\documentclass[12pt,a4paper]{article}
\usepackage{amscd,amsmath,amssymb,amsfonts,times}
\usepackage{mathrsfs}
\usepackage[dvips]{color} 
\usepackage[all]{xy}
\usepackage{enumerate}

\numberwithin{equation}{section}
    \newtheorem{thm}{Theorem}[section]
    \newtheorem{lem}[thm]{Lemma}
    
    \newtheorem{prop}[thm]{Proposition}
    \newtheorem{cor}[thm]{Corollary}

    \newtheorem{conj}[thm]{Conjecture}
    
    \newtheorem{exmp}[thm]{Example}
    \newtheorem{rem}[thm]{Remark}
\DeclareMathAlphabet{\mathpzc}{OT1}{pzc}{m}{it}
\newcommand{\qed}
{\mbox{}\nolinebreak$\square$\medbreak\par}
\newenvironment{pf}{\par\smallskip\noindent\emph{Proof.}}{\hfill\qed\par\smallskip}
\newenvironment{pf*}[1]{\par\smallskip\noindent\emph{#1.}}{\hfill\qed\par\smallskip}
\newcommand{\bysame}{\hskip.3em \leavevmode\rule[.5ex]{2.5em}{.3pt}\hskip0.5em}
\begin{document}
\title{A generalization of the Ross symbols in higher $K$-groups and hypergeometric functions I}\author{M. Asakura
\footnote{
Hokkaido University, 
Sapporo 060-0810, JAPAN. E-mail:
\texttt{asakura@math.sci.hokudai.ac.jp}
}}
\date\empty
\maketitle

\def\can{\mathrm{can}}
\def\cano{\mathrm{canonical}}
\def\ch{{\mathrm{ch}}}
\def\Coker{\mathrm{Coker}}
\def\crys{\mathrm{crys}}
\def\dlog{{\mathrm{dlog}}}
\def\dR{{\mathrm{d\hspace{-0.2pt}R}}}            
\def\et{{\mathrm{\acute{e}t}}}  
\def\Frac{{\mathrm{Frac}}}
\def\phami{\phantom{-}}
\def\id{{\mathrm{id}}}              
\def\Image{{\mathrm{Im}}}        
\def\Hom{{\mathrm{Hom}}}  
\def\Ext{{\mathrm{Ext}}}
\def\MHS{{\mathrm{MHS}}}  
  
\def\can{\mathrm{can}}
\def\arg{\mathrm{arg}}
\def\ch{{\mathrm{ch}}}
\def\Coker{\mathrm{Coker}}
\def\crys{\mathrm{crys}}
\def\dlog{d{\mathrm{log}}}
\def\dR{{\mathrm{d\hspace{-0.2pt}R}}}            
\def\et{{\mathrm{\acute{e}t}}}  
\def\Frac{\operatorname{Frac}}
\def\phami{\phantom{-}}
\def\id{{\mathrm{id}}}              
\def\Image{{\mathrm{Im}}}        
\def\Hom{\operatorname{Hom}}  
\def\Ext{{\mathrm{Ext}}}
\def\MHS{{\mathrm{MHS}}}  
  
\def\ker{\operatorname{Ker}}          
\def\mf{{\text{mapping fiber of}}}
\def\Pic{{\mathrm{Pic}}}
\def\CH{{\mathrm{CH}}}
\def\NS{{\mathrm{NS}}}
\def\NF{{\mathrm{NF}}}
\def\End{\operatorname{End}}
\def\pr{{\mathrm{pr}}}
\def\red{{\mathrm{red}}}
\def\Proj{\operatorname{Proj}}
\def\ord{\operatorname{ord}}
\def\rig{{\mathrm{rig}}}
\def\reg{{\mathrm{reg}}}          %
\def\res{{\mathrm{res}}}          %
\def\Res{\operatorname{Res}}
\def\Spec{\operatorname{Spec}}     
\def\syn{{\mathrm{syn}}}
\def\cont{{\mathrm{cont}}}
\def\ln{{\operatorname{ln}}}

\def\bA{{\mathbb A}}
\def\bC{{\mathbb C}}
\def\C{{\mathbb C}}
\def\G{{\mathbb G}}
\def\bE{{\mathbb E}}
\def\bF{{\mathbb F}}
\def\F{{\mathbb F}}
\def\bG{{\mathbb G}}
\def\bH{{\mathbb H}}
\def\bJ{{\mathbb J}}
\def\bL{{\mathbb L}}
\def\cL{{\mathscr L}}
\def\bN{{\mathbb N}}
\def\bP{{\mathbb P}}
\def\P{{\mathbb P}}
\def\bQ{{\mathbb Q}}
\def\Q{{\mathbb Q}}
\def\bR{{\mathbb R}}
\def\R{{\mathbb R}}
\def\bZ{{\mathbb Z}}
\def\Z{{\mathbb Z}}
\def\cH{{\mathscr H}}
\def\cD{{\mathscr D}}
\def\cE{{\mathscr E}}
\def\cF{{\mathscr F}}
\def\cU{{\mathscr U}}
\def\O{{\mathscr O}}
\def\cR{{\mathscr R}}
\def\cS{{\mathscr S}}
\def\cV{{\mathscr V}}
\def\cX{{\mathscr X}}
\def\cM{{\mathscr M}}
%
\def\ve{\varepsilon}
\def\vG{\varGamma}
\def\vg{\varGamma}
%
%
%
%
\def\lra{\longrightarrow}
\def\lla{\longleftarrow}
\def\Lra{\Longrightarrow}
\def\hra{\hookrightarrow}
\def\lmt{\longmapsto}
\def\ot{\otimes}
\def\op{\oplus}
\def\l{\lambda}
\def\Isoc{{\mathrm{Isoc}}}
\def\Fil{{\mathrm{Fil}}}

\def\MHS{{\mathrm{MHS}}}
\def\MHM{{\mathrm{MHM}}}
\def\VMHS{{\mathrm{VMHS}}}
\def\sm{{\mathrm{sm}}}
\def\tr{{\mathrm{tr}}}
\def\fib{{\mathrm{fib}}}
\def\FIsoc{{F\text{-Isoc}}}
\def\FMIC{{F\text{-MIC}}}
\def\Log{{\mathscr{L}{og}}}
\def\FilFMIC{{\mathrm{Fil}\text{-}F\text{-}\mathrm{MIC}}}

\def\wt#1{\widetilde{#1}}
\def\wh#1{\widehat{#1}}
\def\spt{\sptilde}
\def\ol#1{\overline{#1}}
\def\ul#1{\underline{#1}}
\def\us#1#2{\underset{#1}{#2}}
\def\os#1#2{\overset{#1}{#2}}

\def\s{{(\sigma)}}
\def\Ross{{\xi_{\mathrm{Ross}}}}
\def\HG{{\mathrm{HG}}}


\begin{abstract}
The Ross symbol is defined to be an element $\{1-z,1-w\}$ in 
$K_2$ of a Fermat curve $z^n+w^m=1$.
Ross showed that it is non-torsion by computing the Beilinson regulator.
In this paper,
we introduce a certain generalization of the Ross symbols in $K_{d+1}$ 
of a variety $(1-x_0^{n_0})\cdots(1-x_d^{n_d})=t$.
The main result is that the Beilinson regulator is described by the hypergeometric functions
${}_{d+3}F_{d+2}$'s.
\end{abstract}
\section{Introduction}
In his paper \cite{ross2}, R. Ross introduced an
element 
\begin{equation}\label{intro-eq1}
\{1-z,1-w\}
\end{equation}
in $K_2$ of a Fermat curve $F$ defined by an equation $z^n+w^n=1$.
He proved that it is non-torsion by showing that the real regulator
\begin{equation}\label{intro-eq2}
\reg_\R(\{1-z,1-w\})\in H^2_\cD(F,\R(2))\cong H^1_B(F,\R(1))
\end{equation}
does not vanish in the Deligne-Beilinson cohomology
$H^\bullet_\cD$ where
$H_B^\bullet$ denotes the Betti cohomology.
In this paper we call his element the {\it Ross symbol}.
This is integral in the sense of \cite{Scholl}, and then one can discuss the
regulators and the special values of $L$-functions according to the Beilinson conjecture
\cite{beilinson}, \cite{schneider}. This is studied by Ross \cite{ross1}, K. Kimura \cite{kimurak} and N. Otsubo \cite{otsubo-1}, \cite{otsubo-2}.

\medskip

The purpose of this paper is to introduce a {\it higher Ross symbol} which is
a generalization of the Ross symbol in higher $K$-groups.
Let $A$ be a commutative ring.
Let $n_0,\ldots,n_d$ be positive integers. 
We call an affine scheme 
\[
U=\Spec A[x_0,\ldots,x_d]/((1-x_0^{n_0})\cdots(1-x_d^{n_d})-t),\quad (t\in A)
\]
a {\it hypergeometric scheme} (\S \ref{setting-sect}). 
A particular case is $(1-x_0^{n_0})(1-x_1^{n_1})=1$
$\Leftrightarrow$ $x_0^{-n_0}+x_1^{-n_1}=1$ the Fermat curve.
When $A=\C$,
the periods of integrals are given
by the hypergeometric functions (Theorem \ref{per-thm1}),
\[
{}_{d+1}F_d\left({a_0,\ldots,a_d\atop1,\ldots,1};t\right):=\sum_{n=0}^\infty
\frac{(a_0)_n}{n!}\cdots\frac{(a_d)_n}{n!}t^d,\quad (\alpha)_n:=\alpha(\alpha+1)\cdots(\alpha+n-1)
\]
that is why we call the hypergeometric scheme.
We refer \cite{slater} or \cite[15,16]{NIST}
for the general theory of hypergeometric functions.
For $n_i$-th roots $\nu_i$ of unity,
we define the higher Ross symbol to be a symbol
\[
\Ross=\left\{\frac{1-x_0}{1-\nu_0x_0},\ldots,\frac{1-x_d}{1-\nu_dx_d}\right\}\in K^M_{d+1}(\O(U))
\]
in Milnor's $K$-group. 
See \S \ref{Ross-defn-sect} for the comparison with the original Ross symbol \eqref{intro-eq1}.
We note that, in case $d=1$, this symbol was already discussed in \cite{A}.

\medskip

Otsubo \cite{otsubo-1}, \cite{otsubo-2}
discovered that the real regulators \eqref{intro-eq2} for Fermat curves
can be described in terms of
values at $x=1$ of the hypergeometric functions ${}_3F_2$'s.
A relevant but extended formula is provided in \cite{A},
where the regulators of the higher Ross symbols for the curve
$(1-x_0^{n_0})(1-x_1^{n_1})=t$ are described 
in terms of ${}_4F_3(t)$'s (loc.cit. Thoerem 3.2).
The main result of this paper is a generalization of these formulas for higher dimension.
More precisely, let
\[
F_{a_0,\ldots,a_d}(t):={}_{d+3}F_{d+2}\left({a_0+1,\ldots,a_d+1,1,1\atop 2,\ldots,2};t\right),
\]
be the hypergeometric function, and put
\[
\cF_{a_0,\ldots,a_d}(t):=\sum_{k=0}^d(\psi(a_k)+\gamma)+\log(t)+a_0\cdots a_d\,t\, 
F_{a_0,\ldots,a_d}(t)
\]
where $\gamma=-\Gamma'(1)$ is the Euler constant and
$\psi(t)=\Gamma'(t)/\Gamma(t)$ is the digamma function, cf. \S \ref{analytic-sect}.
Let
$U$ be the hypergeometric scheme
over $\C$ with $t\in \C\setminus\{0,1\}$.
Let
\[\xymatrix{
\reg_B:K^M_{d+1}(\O(U))\ar[r]&H^{d+1}_\cD(U,\Q(d+1))\cong H^d_B(U,\C/\Q(d+1))
}\]
be the symbol map which is a particular case of 
the Beilinson regulator map (e.g. \cite{schneider}).
Then the main theorem is the following.
\begin{thm}[Theorem \ref{mainB-thm}]\label{intro-thm}
Suppose $n_k>1$ for all $k$.
There is a certain homology cycle $\Delta\in H_d(U,\Z)$ such that
the pairing 
$\langle\reg_B(\Ross)\mid\Delta\rangle$ is
\[
\sum_{0< i_k<n_k}(1-\nu_0^{i_0})\cdots(1-\nu_d^{i_d})
\frac{(2\pi i)^d}{n_0\cdots n_d}\cF_{a_0\ldots a_d}(t)\mod\Q(d+1)
\]
where $a_k=1-i_k/n_k$.
\end{thm}
Here we note about the boundary of the higher Ross symbols.
Let $X\supset U$ be a smooth compactification (cf. Proposition \ref{sc-thm}).
We expect that $\Ross$ has no boundary, namely it
lies in the image of Quillen's $K$-group $K_{d+1}(X)^{(d+1)}$
(see \S \ref{boundary-sect}).
This is true in case $d=1,2$ (Corollary \ref{boundary-prop1}), 
while the author has not succeeded to prove it in general.
However it is worth noticing that
one can show that $\reg_B(\Ross)$ lies in the image of
$H^d_B(X,\C/\Q(d+1))$ (see \eqref{mainB-ext-d}),
\[
\reg_B(\Ross)\in \Image[H^d_B(X,\C/\Q(d+1))\to H^d_B(U,\C/\Q(d+1))].
\]

\medskip

Theorem \ref{intro-thm} has some applications to the study of the Beilinson conjecture.
In particular, employing a formula of D. Samart \cite{samart1}, we can prove the following
theorem.
\begin{thm}[Theorem \ref{samart-thm}]
Let $S$ be the K3 surface over $\Q$ 
defined by an equation $(1-x_0^2)(1-x_1^2)(1-x_2^2)=1$.
Let
\[
\Ross=\left\{\frac{1-x_0}{1+x_0},\frac{1-x_1}{1+x_1},\frac{1-x_2}{1+x_2}\right\}
\]
be a higher Ross symbol, lying 
in the image of $K_3(S)$ by Corollary \ref{boundary-prop1}. Then
\[
\frac{1}{(2\pi\sqrt{-1})^2}\langle\reg_\R(\Ross)\mid\Delta^+_1\rangle
=-8L'(h^2_\tr(S),0),
\]
where the notation is as in \S \ref{add-sect}.
\end{thm}
We note that there is a generically finite dominant morphism 
from the self-product of the elliptic curve $x^2+y^4=1$ onto $S$, and 
the $L$-function $L(h^2_\tr(S),s)$ agrees with $L(A,s)$ where $A=\eta^6(4z)$ 
(see \S \ref{add-sect} for the detail). 

\medskip

In the sequel paper \cite{As-Ross2}, 
we give a $p$-adic counterpart of Theorem \ref{intro-thm}
where the $p$-adic analytic function
$\cF^{(\sigma)}_{a_0,\ldots,a_d}(t)$
introduced by the author \cite{New} plays the corresponding role to 
the complex analytic function $\cF_{a_0,\ldots,a_d}(t)$.

\bigskip

The paper is organized as follows.
In \S \ref{HGF-sect}, we introduce hypergeometric schemes. 
We construct a smooth projective compactification over an arbitrary ring $A$
(Proposition \ref{sc-thm}).
In \S \ref{per-sect} we show that periods of integrals are given by
the hypergeometric functions, which is a fundamental formula on hypergeometric schemes. 
\S \ref{dim-sect} is devoted to compute the cohomology groups.
The results in \S \ref{dim-sect} plays a key role in the proof of the main theorem.
The higher Ross symbols are introduced in \S \ref{HRoss-sect}.
In \S \ref{Beilinson-sect} we prove
the main theorem (=Theorem \ref{mainB-thm}) after introducing 
$\cF_{a_0\ldots,a_d}(t)$ and its connection formula (=Theorem \ref{conn})
in \S \ref{analytic-sect}.
Finally we apply the main theorem to the Beilinson conjecture in \S \ref{application-sect}.



\section{Hypergeometric Schemes}\label{HGF-sect}
\subsection{Definition}\label{setting-sect}
Let $d\geq 1$ and 
let $n_i\geq 1$ be positive integers for $i=0,\ldots, d$.
Let $A$ be an integral domain such that all $n_i$ are invertible.
For $t\in A$, let 
\[
U=U_t:=\Spec A[x_0,\ldots,x_d]/((1-x_0^{n_0})\cdots(1-x_d^{n_d})-t)
\] 
be an affine scheme. We call this an {\it hypergeometric scheme} over $A$.
If $t(1-t)\in A^\times$, then $U$ is smooth over $A$.

We denote by $\mu_{m}=\mu_{m}(A)\subset A^\times$ 
the group of $m$-th roots of unity in $A$.
Let $G=\mu_{n_0}\times \cdots\times \mu_{n_d}$.
We write
$\sigma_j(\nu_j):=(1,\ldots,\nu_j,\ldots,1)\in G$ and 
$\sigma(\ul\nu):=\sigma_0(\nu_0)
\cdots\sigma_d(\nu_d)$
for $\ul \nu=(\nu_0,\ldots,\nu_d)\in G$.
The group $G$ acts on 
$U$ in a way that
\begin{equation}\label{setting-eq1}
\sigma(\ul \nu)(x_0,\ldots,x_d)=(\nu_0x_0,\ldots,\nu_dx_d).
\end{equation} 
For a commutative ring $k$ and
a $k[G]$-module $H$ and a homomorphism
$\chi:G\to k^\times$,
let
\[
H(\chi):=\{x\in H\mid \sigma(\ul\nu)(x)=
\chi(\ul\nu)x,\,\forall\, \ul\nu=(\nu_0,\ldots,\nu_d)
\in G\}
\]
denote the eigenspace where $\chi(\ul\nu)=\chi(\nu_0)
\cdots\chi(\nu_d)$. If $A$ is an integral $k$-algebra and $k$ is a field which
contains primitive $n_i$-th roots of unity for each $i$, then
there is a natural correspondence
\[
\Z/n_0\Z\times\cdots\times\Z/n_d\Z\os{\cong}{\lra}
\Hom(G,k^\times),\quad
(i_0,\ldots,i_d)\longmapsto\chi(i_0,\ldots,i_d)
\]
where $\chi(i_0,\ldots,i_d)$ is given by 
$\chi(i_0,\ldots,i_d)(\nu_0,\ldots,\nu_d)=\nu_0^{i_0}\cdots\nu_d^{i_d}$.
We also write $H(i_0,\ldots,i_d)=H(\chi(i_0,\ldots,i_d))$ simply.

\subsection{Smooth Compactification}
We construct a smooth compactification of $U/A$.
Let $\P^1_A(x_i)=\mathrm{Proj}A[X_i,Y_i]$ denote the projective line over $A$
with inhomogeneous coordinate $x_i=X_i/Y_i$.
Let
\begin{equation}\label{sc-eq1}
X^*=X^*_t\subset 
\P^1_A(x_0)\times \cdots\times \P^1_A(x_d)
\end{equation}
be a closed subscheme defined by an equation
\[
(Y_0^{n_0}-X_0^{n_0})\cdots(Y_d^{n_d}-X_d^{n_d})=tY_0^{n_0}\cdots Y_d^{n_d}.
\]
There is an open embedding $U\hra X^*$.

\begin{prop}\label{sc-thm}
Suppose that $t(1-t)\in A^\times$.
Then there is a projective morphism
\[
\rho:X\lra X^*\subset \P^1_A(x_0)\times \cdots\times \P^1_A(x_d)
\]
such that $X\to \Spec A$ is smooth, 
$\rho^{-1}(U)\os{\sim}{\to} U$ and
$Z:=X\setminus\rho^{-1}(U)$ is a relative simple NCD over $A$.
Moreover any geometric fiber of $X\to\Spec A$ is connected.
\end{prop}
There is {\it no} natural smooth compactification of $U$.
The author does not know whether the action of $G$ extends or not.
\begin{pf}
The singular locus (= the closed subset which is not smooth over $A$) of $X^*$
lies outside $U$.
Put $X^0:=X^*$ and $Z^0:=X^*\setminus U$.
A neighborhood of $Z^0$ in $X^0$ is locally described by an equation 
\[
(\nu_{i_1}-x_{i_1})\cdots(\nu_{i_r}-x_{i_r})=\text{(unit)}\times y_{j_1}^{n_{j_1}}\cdots
y_{j_s}^{n_{j_s}}
\]
and the singular locus is the union of 
\[
\{\nu_i-x_i=\nu_j-x_j=y_k=0\}\cong (\P^1_A)^{d-2}.  
\]
We denote the above situation by
\[
{\mathbb A}^{d+1}(z,w)\supset 
\overbrace{\{z_1\cdots z_r=uw_1^{m_1}\cdots w_s^{m_s}\}}^{\text{locus of }X^0}
\supset
\overbrace{\{w_1\cdots w_s=0\}}^{\text{locus of }Z^0}.
\]
We construct a sequence of blow-ups
\begin{equation}\label{blow-up-eq0}
\cdots \lra X^{i+1}\lra X^i\lra\cdots\lra X^0=X^*
\end{equation}
in the following way. Let $Z^i\subset X^i$ be the reduced inverse image $Z^0$.
Choose an irreducible component $Z_0^i$ of $Z^i$ which contains at least one component
of the singular locus of $X^i$. Then we take the blow-up $X^{i+1}\to X^i$ along
$Z^i_0$.
We claim that each $X^i\supset Z^i$ 
is locally (in Zariski topology) described as
\begin{equation}\label{blow-up-eq1}
{\mathbb A}^{d+1}(z,w,v)\supset 
\overbrace{\{z_1\cdots z_r=uw_1^{m_1}\cdots w_s^{m_s}\}}^{\text{locus of }X^i}
\supset
\overbrace{\{w_1\cdots w_sv_1\cdots v_l=0\}}^{\text{locus of }Z^i}
\end{equation}
with $u$ a unit (possibly $l=0$). 
If $i=0$, this is straightforward.
Suppose that $X^i\supset Z^i$ is described as above.
Put $Z^i_{ab}:=\{z_a=w_b=0\}$ and $Z^i_c:=\{v_c=0\}$
irreducible components of $Z^i$.
If $X^{i+1}\to X^i$ be the blow-up along $Z^i_c$, then this is an isomorphism
over the above locus, so there is nothing to prove.
Suppose that $\rho:X^{i+1}\to X^i$ is the blow-up along $Z^i_{ab}$.
We may assume $a=b=1$ without loss of generality.
Then $X^{i+1}$ is covered by two affine open sets ($r+s+l=d+1$)
\begin{align}
U_1=&\left\{
\left(z_1,\ldots,z_r,\frac{w_1}{z_1},w_2,\ldots,w_s,v_1,\ldots,v_l\right)\bigg|
z_2\cdots z_r=u_1z_1^{m_1-1}\left(\frac{w_1}{z_1}\right)^{m_1}
w_2^{m_2}\cdots w_s^{m_s}
\right\}\label{blow-up-eq2}
\\
U_2=&\left\{\left(\frac{z_1}{w_1},z_2,\ldots,z_r,w_1,\ldots,w_s,v_1,\ldots,v_l\right)\bigg|
\frac{z_1}{w_1}\cdot z_2\cdots z_r=u_2
w_1^{m_1-1}w_2^{m_2}\cdots w_s^{m_s}
\right\}\label{blow-up-eq3}
\end{align}
with $u_i$ units, and
\[
\rho^{-1}(Z^i)\cap U_1=
\left\{
z_1\cdot\frac{w_1}{z_1}\cdot w_2\cdots w_s=0\right\},\quad
\rho^{-1}(Z^i)\cap U_2=
\left\{
 w_1\cdots w_s=0\right\}.
\]
Hence $X^{i+1}\supset Z^{i+1}$ is also described as in \eqref{blow-up-eq1}.
In the sequence \eqref{blow-up-eq0}, $r$ and $\max\{m_i\}$ are decreasing sequences.
Thus $X^i\supset Z^i$ for $i\gg0$ is locally described by either of the following descriptions
\begin{enumerate}
\item[(a)]
$\bA^{d+1}(z,w,v)\supset\{z_1\cdots z_r=uw_1\cdots w_s\}\supset
\{w_1\cdots w_sv_1\cdots v_l=0\}$, $r,s\geq2$
\item[(b)]
$\bA^{d+1}(z,w,v)\supset\{z_1\cdots z_r=uw_1\}\supset\{w_1v_1\cdots v_l=0\}$,
\item[(c)]
$\bA^{d+1}(z,w,v)\supset\{z_1=uw^{m_1}_1\cdots w^{m_s}_s\}\supset
\{w_1\cdots w_sv_1\cdots v_l=0\}$,
\item[(d)]
$\bA^{d+1}(z,w,v)\supset\{z_1\cdots z_r=u\}\supset
\{w_1\cdots w_sv_1\cdots v_l=0\}$.
\end{enumerate}
If a local description of $X^i\supset Z^i$ is either of (b), (c) or (d), 
then $X^i$ is smooth over $A$ and $Z^i$ is a relative NCD.
Moorever if we take the blow-up $\rho:X^{i+1}\to X^i$ along a component of $Z^i$, it is an
isomorphism over the above locus as $\mathrm{codim}(Z^i)=1$,
so the local description remains the same.
Consider that a local description is as in case (a).
We keep continuing the blow-ups. Then $r+s$ is strictly decreasing, so that 
the description will finally be the case (b), (c) ($m_i=1$) or (d).
We thus have a sequence of blow-ups
\[
X=X^N\lra\cdots \lra X^{i+1}\os{\rho_{i+1}}{\lra}X^i\to\cdots \lra X^0=X^*
\]
such that the local descriptions of $X\supset Z$ (in Zariski topology)
are the case (b), (c) or (d).
By the construction, each $X^i$ has an affine covering $\cup V^i_\alpha$ 
such that the geometric fibers of $\rho_{i+1}^{-1}(V^i_\alpha)\to\Spec A$ are irreducible.
Since the geometric fibers of of $X^*\to\Spec A$ are irreducible, 
the same thing holds for every $X^i$.
In particular $X\to\Spec A$ is a projective smooth morphism with connected fibers.
Finally we see that each irreducible component of $Z$ is smooth over $A$.
Indeed the local descriptions (b), (c) or (d) tells that
there is an affine covering $X=\cup_\alpha V_\alpha$ 
such that each $Z\cap V_\alpha$ is a relative simple NCD over $A$.
Let $Z_j$ be an arbitrary irreducible component of $Z$.
Since $Z_j\cap V_\alpha$
is irreducible (notice that any non-empty Zariski open set of an irreducible scheme is irreducible), 
it turns out that $Z_j\cap V_\alpha$ is smooth over $A$
for any $\alpha$, which means that $Z_j$ is smooth
over $A$.
Hence $Z$ is a relative simple NCD over $A$.
This completes the proof.
\end{pf}

\subsection{Note on Boundary components}\label{tate-sect}
Let $V$ be a quasi-projective $A$-scheme.
Let $l$ be a prime invertible in $A$.
Define the etale homology $H_j^\et(V,\Q_l)$
to be the sheaf 
$\Q\ot\varprojlim_n R^{2n-j}\pi_*i^!\Z/l^n$ on $\Spec A$ where 
$i:V\hra P$ is a closed immersion to a smooth $A$-scheme $P$ of relative dimension $n$
and $\pi:V\to\Spec A$ is the structural morphism.
We call $V$ a {\it mixed Tate motive} if
$H_k^\et(V,\Q_l)$ 
are succesive extensions of products of $\Q_l(j)$'s for any $k\in \Z_{\geq 0}$ and 
primes $l$ invertible in $A$.
A mixed Tate motive $V$ is called a {\it Tate motive} if $\pi$ is projective and smooth.

\begin{prop}\label{tate-prop}
Let $X$ be the smooth compactification in Proposition \ref{sc-thm} and 
$Z=X\setminus U$ the boundary.
Let $Z=\cup_k Z_k$ be the decomposition where $Z_k$ are $A$-smooth
divisors.
Then any intersection $Z_{k_1}\cap\cdots\cap Z_{k_q}$ is a Tate motive over $A$.
\end{prop}
\begin{pf}
Recall from the proof of Proposition \ref{sc-thm} the sequence of the blow-ups
$\cdots\to X^{i+1}\to X^i\to\cdots$.
We use the notation $Z^i=\bigcup Z^i_{ab}\cup\bigcup Z^i_c$.
For a finite set $I$ of indices of $ab$ and $c$'s, write
\[
(Z^i)^I=Z^i_{a_1b_1}\cap\cdots Z^i_{a_nb_n}\cap Z^i_{c_1}\cap\cdots\cap Z^i_{c_m}.
\]
It is straightforward to see that 
every intersections $(Z^0)^I$ are a mixed Tate motive.
We show that this is true for all $i$ by the induction.
To do this, it is enough to show that
$(Z^{i+1})^I$ satisfies either of the following.
Let $\rho:X^{i+1}\to X^i$ be the blow-up.
\begin{itemize}
\item[(Z1)] $(Z^{i+1})^I=\emptyset$,
\item[(Z2)]
$(Z^{i+1})^I\os{\cong}{\to} (Z^i)^J$ with some $J$,
\item[(Z3)] $(Z^{i+1})^I$ is a $\P^1$-bundle over $(Z^i)^J$ with some $J$,
\item[(Z4)] 
there are $(Z^i)^J$ and a closed subscheme $E\subset
(Z^i)^J$ which is a union of $(Z^i)^{J'}$'s
such that $\rho((Z^{i+1})^I)=(Z^i)^J$ and 
$(Z^{i+1})^I\setminus \rho^{-1}(E)\os{\cong}{\to} (Z^i)^J\setminus E$ and
$\rho^{-1}(E)\to E$ is a $\P^1$-bundle.
\end{itemize}
If the above holds, then it follows that 
$(Z^{i+1})^I$ is a mixed Tate motives which can be shown with use of
the localization sequence
\[\xymatrix{
\cdots\ar[r]&H_i^\et(D,\Q_l)\ar[r]&H_i^\et(V,\Q_l)\ar[r]&H_i^\et(V\setminus D,\Q_l)\ar[r]&
\cdots
}\]
of etale homology groups
for quasi-projective $A$-schemes $(D,V)$ with $D\subset V$ a closed subscheme.

\medskip

The blow-up $\rho:X^{i+1}\to X^i$ is locally described by either of the following cases.
\begin{enumerate}
\item[(i)]
$\rho:X^{i+1}\to X^i$ is the blow-up along $Z^i_{ab}=\{z_a=w_b=0\}$, and 
$m_a>1$,
\item[(ii)]
$\rho:X^{i+1}\to X^i$ is the blow-up along $Z^i_{ab}=\{z_a=w_b=0\}$, and 
$m_a=1$,
\item[(iii)]
$\rho:X^{i+1}\to X^i$ is the blow-up along $Z^i_{c}=\{v_c=0\}$.
\end{enumerate}
In the case (iii), 
$\rho$ is an isomorphism over the locus \eqref{blow-up-eq1}, and hence
$Z^{i+1}\cong Z^i$.
Consider the case (i).
For simplicity of the notation, we take $a=b=1$.
Let $X^{i+1}=U_1\cup U_2$ be the affine covering in 
\eqref{blow-up-eq2} and \eqref{blow-up-eq3}.
If $m_1\geq 2$ then
\begin{align*}
\rho^{-1}(Z^i_{11})\cap U_1&=
\left\{
\left(z_1,\ldots,z_r,\frac{w_1}{z_1},w_2,\ldots,w_s,v_1,\ldots,v_l\right)\bigg|
z_1=z_2\cdots z_r=0
\right\}\\
\rho^{-1}(Z^i_{11})\cap U_2&=
\left\{\left(\frac{z_1}{w_1},z_2,\ldots,z_r,w_1,\ldots,w_s,v_1,\ldots,v_l\right)\bigg|
w_1=\frac{z_1}{w_1}\cdot z_2\cdots z_r=0\right\}.
\end{align*}
The exceptional divisors are $E^{i+1}_a$, $2\leq a\leq r$ where
\[
E^{i+1}_a\cap U_1=\{z_1=z_a=0\},\quad
E^{i+1}_a\cap U_2=\{z_a=w_1=0\}
\]
is a $\P^1$-bundle over $Z^i_{11}\cap Z^i_{a1}$.
The proper transform $Z^{i+1}_{1b}$ of $Z^i_{1b}$ with $b\geq 1$ is given by
\[
Z^{i+1}_{1b}\cap U_1=\emptyset,\quad
Z^{i+1}_{1b}\cap U_2=\{\frac{z_1}{w_1}=w_b=0\},
\]
the proper transform $Z^{i+1}_{a1}$ of $Z^i_{a1}$ with $a\geq 2$ is given by
\[
Z^{i+1}_{a1}\cap U_1=\{z_a=\frac{w_1}{z_1}=0\},\quad
Z^{i+1}_{a1}\cap U_2=\emptyset.
\]
Hence $Z^{i+1}_{ab}\cong Z^{i}_{ab}$ if $a=1$ or $b=1$.
The proper transform $Z^{i+1}_{ab}$ of $Z^i_{ab}$ with $a,b\geq 2$ is given by
\[
Z^{i+1}_{ab}\cap U_1=\{z_a=w_b=0\},\quad
Z^{i+1}_{ab}\cap U_2=\{z_a=w_b=0\}.
\]
Hence 
$Z^{i+1}_{ab}\to Z^{i}_{ab}$ is the blow-up along $Z^i_{ab}\cap Z^i_{11}$.
One easily sees that the proper transform $Z^{i+1}_{c}$ of $Z^i_{c}$ is also isomorphic to
$Z^i_c$.
\[
Z^{i+1}:=\rho^{-1}(Z^i)=\bigcup_{a\geq 2}E^{i+1}_a
\cup\bigcup_{b\geq 1}Z^{i+1}_{1b}\cup\bigcup_{a\geq 2}Z^{i+1}_{a1}
\cup\bigcup_{a,b\geq2}Z^{i+1}_{ab}.
\]
By the above descriptions, we see
\[
Z^{i+1}_{1b}\cap Z^{i+1}_{1b'}\cong
Z^{i}_{1b}\cap Z^{i}_{1b'},\quad
Z^{i+1}_{1b}\cap Z^{i+1}_{a1}=\emptyset,\quad
Z^{i+1}_{a1}\cap Z^{i+1}_{a'1}\cong
Z^{i}_{a1}\cap Z^{i}_{a'1},
\]
for $b,b'\geq 1$ and $a,a'\geq 2$,
\[
Z^{i+1}_{1b}\cap Z^{i+1}_{a'b'}\cong
Z^{i}_{1b}\cap Z^{i}_{a'b'},\quad
Z^{i+1}_{a1}\cap Z^{i+1}_{a'b'}\cong
Z^{i}_{a1}\cap Z^{i}_{a'b'},
\]
for $b\geq 1$ and $a,a',b'\geq 2$,
and 
$Z^{i+1}_{ab}\cap Z^{i+1}_{a'b'}\to Z^{i}_{ab}\cap Z^{i}_{a'b'}$ is the blow-up along 
$Z^i_{ab}\cap Z^{i}_{a'b'}\cap Z^i_{11}$ for $a,a',b,b'\geq 2$.
Moreover
\[
E^{i+1}_a\cap Z^{i+1}_{1b}\cong
Z^i_{a1}\cap Z^{i}_{1b},\quad
E^{i+1}_a\cap Z^{i+1}_{a'1}\cong
Z^i_{11}\cap Z^{i}_{a1}\cap Z^{i}_{a'1},\quad
\]
for $a,a'\geq 2$ and $b\geq 1$, and 
$E^{i+1}_{a}\cap Z^{i+1}_{a'b'}\to Z^i_{11}\cap Z^i_{a1}\cap Z^{i}_{a'b'}$ is a $\P^1$-bundle 
for $a,a',b'\geq 2$.
Besides
$Z^{i+1}_c\cap Z^{i+1}_{ab}\cong Z^{i}_c\cap Z^{i}_{ab}$
if $a=1$ or $b= 1$, and
$Z^{i+1}_c\cap Z^{i+1}_{ab}\to Z^{i}_c\cap  Z^{i}_{ab}$ is the blow-up along 
$Z^{i}_c\cap Z^i_{ab}\cap Z^i_{11}$ if $a,b\geq 2$.
In all cases, the descriptions (Z1),\ldots,(Z4) remain true.

\medskip

Next we consider the case (ii), namely
$\rho:X^{i+1}\to X^i$ is the blow-up along $Z^i_{11}$ and $m_1=1$.
Then
\begin{align*}
\rho^{-1}(Z^i_{11})\cap U_1&=
\{z_1=0,\, z_2\cdots z_r=u_0\cdot\frac{w_1}{z_1}\cdot
w_2^{m_2}\cdots w_s^{m_s}\},\\
\rho^{-1}(Z^i_{11})\cap U_2&=
\{w_1=0,\, \frac{z_1}{w_1}\cdot z_2\cdots z_r=u_0\cdot
w_2^{m_2}\cdots w_s^{m_s}\},
\end{align*}
where $u_0:=u|_{z_1=w_1=0}$.
This shows that $\rho^{-1}(Z^i_{11})\to Z^i_{11}$ is the blow-up along the ideal
$(z_2\cdots z_r,w_2^{m_2}\cdots w_s^{m_s})$.
Hence $\rho^{-1}(Z^i_{11})$ is irreducible (possibly non-smooth over $A$), and letting
$E^i_{11}:=\cup_{a,b\geq 2}Z^i_{11}\cap Z^i_{ab}$, we have that
$\rho^{-1}(Z^i_{11}\setminus E_{11}^i)
\os{\sim}{\to} Z^i_{11}\setminus E_{11}^i$
and $\rho^{-1}(E^i_{11})\to E^i_{11}$
is a $\P^1$-bundle.
There is no exceptional divisor.
The proper transform $Z^{i+1}_{1b}$ of $Z^i_{1b}$ for $b\geq 2$ is given by
\[
Z^{i+1}_{1b}\cap U_1=\emptyset,\quad
Z^{i+1}_{1b}\cap U_2=\{\frac{z_1}{w_1}=w_b=0\},
\]
the proper transform $Z^{i+1}_{a1}$ of $Z^i_{a1}$ for $a\geq 2$ is given by
\[
Z^{i+1}_{a1}\cap U_1=\{z_a=\frac{w_1}{z_1}=0\},\quad
Z^{i+1}_{a1}\cap U_2=\emptyset,
\]
and 
the proper transform $Z^{i+1}_{ab}$ of $Z^i_{ab}$ for $a,b\geq 2$ is given by
\[
Z^{i+1}_{ab}\cap U_1=\{z_a=w_b=0\},\quad
Z^{i+1}_{ab}\cap U_2=\{z_a=w_b=0\}.
\]
Therefore 
$Z^{i+1}_{ab}\to Z^i_{ab}$ is the blow-up along $Z^i_{11}\cap Z^i_{ab}$ if $a,b\geq 2$,
and otherwise $Z^{i+1}_{ab}\cong Z^i_{ab}$ unless $(a,b)=(1,1)$.
We have
\[
Z^{i+1}_{1b}\cap Z^{i+1}_{a'b'}\cong
Z^{i}_{1b}\cap Z^{i}_{a'b'},\quad
Z^{i+1}_{a1}\cap Z^{i+1}_{a'b'}\cong
Z^{i}_{a1}\cap Z^{i}_{a'b'},
\]
and $Z^{i+1}_{ab}\cap Z^{i+1}_{a'b'}\to Z^{i}_{ab}\cap Z^{i}_{a'b'}$ is the blow-up
along  $Z^i_{11}\cap Z^{i}_{ab}\cap Z^{i}_{a'b'}$.
Moreover
\[
Z^{i+1}_{11}\cap Z^{i+1}_{1b}=\{\frac{z_1}{w_1}=w_1=w_b=0\}\cong Z^i_{11}\cap Z^i_{1b},
\]
\[
Z^{i+1}_{11}\cap Z^{i+1}_{a1}=\{z_1=z_a=\frac{w_1}{z_1}=0\}\cong Z^i_{11}\cap Z^i_{a1},
\]
and $Z^{i+1}_{11}\cap Z^{i+1}_{ab}\to Z^{i}_{11}\cap Z^{i}_{ab}$ is a $\P^1$-bundle,
$Z^{i+1}_{11}\cap Z^{i+1}_{c}\to Z^{i}_{11}\cap Z^{i}_{c}$ is the blow-up
along the ideal $(z_2\cdots z_r,w_2^{m_2}\cdots w_s^{m_s})$.
Thus the descriptions (Z1),\ldots,(Z4) remain true also in the case (ii).
\end{pf}

\subsection{Rational differential forms in de Rham cohomology}\label{form-sect}
We assume that each $n_k>1$.
For $(i_0,\ldots,i_d)\in \Z^{d+1}$ with $0<i_k<n_k$, put
\begin{align}
\omega_{i_0\ldots i_d}&:=n_0^{-1}x_0^{i_0-n_0}x_1^{i_1-1}\cdots x_d^{i_d-1}
\frac{dx_1\cdots dx_d}{(1-x_1^{n_1})\cdots(1-x_d^{n_d})}\label{form-eq1}\\
&=(-1)^kn_k^{-1}x_0^{i_0-1}\cdots x_k^{i_k-n_k}\cdots x_d^{i_d-1}
\frac{dx_0\cdots \widehat{dx_k}\cdots dx_d}
{(1-x_0^{n_0})\cdots\widehat{(1-x_k^{n_k})}\cdots(1-x_d^{n_d})}\label{form-eq1-1}\\
&\in \vg(U,\Omega^d_{U/A})\notag
\end{align}
a rational $d$-form.
Moreover we put
\begin{equation}\label{form-eq2}
\omega^{(r)}_{i_0\ldots i_d}
:=\left(\frac{1-x_0^{n_0}}{x_0^{n_0}}\right)^r\omega_{i_0\ldots i_d}\in H^d_\dR(U'/A),
\quad U':=U\setminus\{x_0=0\}
\end{equation}
for an integer $r\geq0$.
Let $X\to\Spec A$ be a (fixed) smooth compactification constructed in Proposition \ref{sc-thm}.
\begin{lem}\label{form-lem1}
$\omega_{i_0\ldots i_d}
\in \vg(X,\Omega^d_{X/A})$.
\end{lem}
If $A$ is a $\Q$-algebra, there is the exact sequence
\[
\xymatrix{
\cdots\ar[r]&H^{i-1}_\dR(Z/A)\ar[r]&H^{i}_{c,\dR}(U/A)\ar[r]
&H^{i}_\dR(X/A)\ar[r]&H^{i}_\dR(Z/A)\ar[r]&\cdots.
}
\]
This induces an isomorphism
\[
F^d\cap H^d_{c,\dR}(U/A)\os{\cong}{\lra} F^d\cap H^d_{\dR}(X/A).
\]
where $F^\bullet$ denotes the Hodge filtrartion.
Hence $\omega_{i_0\ldots i_d}$ also defines a cohomology class in $H^d_{c,\dR}(U)$
\begin{pf}
We may assume $A=\Z[1/N][t]$ with $t$ an indeterminate and $N:=n_0\cdots n_d$,
since $\omega_{i_0\ldots i_d}$ is defined over $\Z[1/N][t]$.
Recall the projective scheme $X^*$ in \eqref{sc-eq1}. 
Let $\rho:X\to X^*$ be the desingularization.
Let $y_i=1/x_i$.
A locus of the boundary $X^*\setminus U$ is described by
\[
(\nu_{k_1}-x_{k_1})\cdots(\nu_{k_r}-x_{k_r})=\text{(unit)}\times y_{j_1}^{n_{j_1}}\cdots
y_{j_s}^{n_{j_s}}
\]
and then one can describe $\omega_{i_0\ldots i_d}$ as follows
\[
\text{(regular function)}\times
dy_{j_1}\cdots dy_{j_s}\wedge
\prod_{k\ne j_0,j_1,\ldots, j_s}\frac{dx_k}{\nu_k-x_k}
\]
where $j_0\in\{0,1,\ldots,d\}\setminus
\{j_1,\ldots,j_s\}$ is a fixed integer.
Therefore its pull-back by $\rho$ belongs to the subspace
\[
\rho^*(y_{j_1}\cdots y_{j_s})\cdot
\Omega^d_{X/A}(\log Z), \quad Z:=X\setminus U
\]
in the locus.
Noticing that $\Omega^d_{X/A}(\log Z)=\Omega^d_{X/A}\ot\O( Z)$ is an invertible sheaf
and $\rho^*(y_{j_1}\cdots y_{j_s})\in \O(-Z)$ locally,
this implies that
\[
\rho^*\omega_{i_0\ldots i_d}\in \vg(X\setminus Z^{[2]},\Omega^d_{X/A})
\]
where $Z^{[r]}=\cup_{i_1<\cdots<i_r}Z_{i_1}\cap\cdots\cap Z_{i_r}$.
Since $X$ is regular (as $A=\Z[1/N][t]$) and $\mathrm{codim}_X(Z^{[2]})= 2$, 
one has $\vg(X\setminus Z^{[2]},\Omega^d_{X/A})=\vg(X,\Omega^d_{X/A})$, and hence
the lemma follows.
\end{pf}
\begin{lem}\label{form-lem2}
Suppose that $A$ is a $\Q$-algebra. 
Then
$\omega_{i_0\ldots i_d}^{(r)}
\in \Image[H^d_\dR(X/A)\to H^d_\dR(U'/A)]$.
\end{lem}
\begin{pf}
We may assume $A=\ol\Q[t,(t-t^2)^{-1}]$.
Let $\partial=\frac{d}{dt}$ be the differential operator on the $A$-module $H^d_\dR(U'/A)$
induced from the Gauss-Manin connection, namely $\partial$ is the composition
\begin{equation}\label{form-lem2-eq1}
H^d_\dR(U'/A)\os{\nabla}{\lra}\Omega^1_{A/\ol\Q}\ot_A H^d_\dR(U'/A)\os{\cong}{\lra}
H^d_\dR(U'/A)
\end{equation}
of arrows where the second one is given by $dt\ot x\mapsto x$. Let $D:=t\partial$.
Put
\[
\Omega:=x_1^{i_1-1}\cdots x_d^{i_d-1}
\frac{dx_1\cdots dx_d}{(1-x_1^{n_1})\cdots(1-x_d^{n_d})}
=\bigwedge_{1\leq k\leq d}
\left(\sum_{\zeta_k\in \mu_{n_k}} n_k^{-1}\zeta_k^{i_k}\frac{dx_k}{x_k-\zeta_k}\right)
\in H^d_\dR(U/A).
\]
Since this lies in the image of $(\O(U)^\times)^{\ot d}$ by the dlog-map,
this is annihilated by the differential operator $D$.
Using the equalities
\[dx_0\wedge\Omega
=\left(\frac{1-x_0^{n_0}}{-n_0x_0^{n_0-1}}\right)\frac{dt}{t}\wedge\Omega,\quad
dx_0\wedge\omega_{i_0\ldots i_d}
=\left(\frac{1-x_0^{n_0}}{-n_0x_0^{n_0-1}}\right)\frac{dt}{t}\wedge\omega_{i_0\ldots i_d}
\]
one has
\begin{align*}
D(\omega_{i_0\ldots i_d})
&=D(n_0^{-1}x^{i_0-n_0}\Omega)\\
&=n_0^{-1}(i_0-n_0)x^{i_0-n_0-1}\left(\frac{1-x_0^{n_0}}{-n_0x_0^{n_0-1}}\right)\Omega\\
&=a_0\left(\frac{1-x_0^{n_0}}{x_0^{n_0}}\right)\omega_{i_0\ldots i_d}\quad (a_0:=1-i_0/n_0)\\
&=a_0\omega^{(1)}_{i_0\ldots i_d},
\end{align*}
and
\begin{align}
D(\omega^{(r)}_{i_0\ldots i_d})
&=
D\left(\left(\frac{1-x_0^{n_0}}{x_0^{n_0}}\right)^r\omega_{i_0\ldots i_d}\right)\notag\\
&=r\left(\frac{1-x_0^{n_0}}{x_0^{n_0}}\right)^{r-1}\frac{-n_0}{x_0^{n_0+1}}
\left(\frac{1-x_0^{n_0}}{-n_0x_0^{n_0-1}}\right)\omega_{i_0\ldots i_d}+
\left(\frac{1-x_0^{n_0}}{x_0^{n_0}}\right)^rD(\omega_{i_0\ldots i_d})\notag\\
&=r\omega^{(r)}_{i_0\ldots i_d}+(a_0+r)\omega^{(r+1)}_{i_0\ldots i_d}.\label{form-lem2-eq2}
\end{align}
Note $0<a_0<1$. Hence this implies 
\[
\langle \omega^{(0)}_{i_0\ldots i_d},\ldots,\omega^{(r)}_{i_0\ldots i_d}\rangle_\Q
=\langle \omega_{i_0\ldots i_d},D(\omega_{i_0\ldots i_d}),\ldots,
D^r(\omega_{i_0\ldots i_d})\rangle_\Q\subset H^d_\dR(U'/A).
\]
Since $\omega_{i_0\ldots i_d}\in \vg(X,\Omega^d_{X/A})$ by Lemma \ref{form-lem2}, 
the right hand side lies in the image of $H^d_\dR(X/A)$.
\end{pf}
\begin{lem}\label{form-lem3}
Suppose that $A$ is a $\Q$-algebra. Let $U'=U\setminus\{x_0=0\}$.
Let $\partial:=\frac{d}{dt}$ be the differential operator defined by 
the composition \eqref{form-lem2-eq1}.
Then
\[
t^r\partial^r(\omega_{i_0\ldots i_d})|_{U'}=a_0(a_0+1)\cdots(a_0+r-1)
\omega_{i_0\ldots i_d}^{(r)}\in H^d_\dR(U'/A),\quad r\geq 1.
\]
According to this, we define a lifting 
\[
\wt\omega_{i_0\ldots i_d}^{(r)}:=
\frac{1}{a_0(a_0+1)\cdots(a_0+r-1)}t^r\partial^r(\omega_{i_0\ldots i_d})
\in H^d_\dR(X/A).
\]
\end{lem}
\begin{pf}
Again we may assume $A=\ol\Q[t,(t-t^2)^{-1}]$.
It follows from \eqref{form-lem2-eq2} that one has
\[
(D-r)(\omega^{(r)}_{i_0\ldots i_d})=(a_0+r)\omega^{(r+1)}_{i_0\ldots i_d}
\]
for all $r\geq 0$ and hence
\[
D(D-1)\cdots(D-r+1)(\omega_{i_0\ldots i_d})=
a_0(a_0+1)\cdots(a_0+r-1)\omega_{i_0\ldots i_d}^{(r)}
\]
by the induction on $r$. Now the assertion follows from
an equality $D(D-1)\cdots(D-r+1)=t^r\partial^r$.
\end{pf}
\section{Periods of Hypergeometric Schemes}
In this section we work over the base field $\C$.  Put $A:=\C[t,(t-t^2)^{-1}]$,
and $S:=\Spec A=\P^1(t)\setminus\{t=0,1,\infty\}$.
Let 
\[
U=\Spec A[x_0,\ldots,x_d]/((1-x_0^{n_0})\cdots(1-x_d^{n_d})-t)
\] be the hypergeometric scheme over $A$ defined in \S \ref{setting-sect}.
We think of $U\to S$ to be a fibration of complex manifolds.
For $\alpha\in \C\setminus\{0,1\}$,
we denote by
$U_\alpha$ the fiber at $t=\alpha$.

\subsection{Periods of integrals for Hypergeometric schemes}\label{per-sect}
\begin{lem}\label{lem-per-thm1-eq2}
Let $n\geq0$ be an interger, and $c_1,c_2\in\C$ with $c_2\ne0$.
Then
\begin{equation}\label{per-thm1-eq2}
\frac{1}{2\pi i}\oint_{|x-1|=\ve}\frac{x^{c_1-1}}{(1-x^{c_2})^{n+1}}dx
=-c_2^{-1}\frac{(1-c_1/c_2)_n}{n!},\quad (0<\ve\ll1)
\end{equation}
where $x^{c}$ takes the branch such that $|x^c-1|\ll1$.
\end{lem}
\begin{pf}
We show the lemma by the induction on $n$.
The case $n=0$ is simple.
Let $k\geq 1$ and
suppose that \eqref{per-thm1-eq2} holds for $n=k-1$.
Then
\begin{align*}
\oint_{|x-1|=\ve}\frac{x^{c_1-1}}{(1-x^{c_2})^{k+1}}dx
&=\oint_{|x-1|=\ve}d\left(\frac{x^{c_1-c_2}}{kc_2(1-x^{c_2})^{k}}\right)
+\frac{1-c_1/c_2}{k}\frac{x^{c_1-c_2-1}}{(1-x^{c_2})^k}dx\\
&=\frac{1-c_1/c_2}{k}\oint_{|x-1|=\ve}\frac{x^{c_1-c_2-1}}{(1-x^{c_2})^k}dx\\
\end{align*}
so that we have
\[
\frac{1}{2\pi i}\oint_{|x-1|=\ve}\frac{x^{c_1-1}}{(1-x^{c_2})^{k+1}}dx
=\frac{1-c_1/c_2}{k}\cdot(-c_2^{-1})\frac{(1-(c_1-c_2)/c_2)_{k-1}}{(k-1)!}
=-c_2^{-1}\frac{(1-c_1/c_2)_k}{k!}
\]
as required.
\end{pf}
\begin{thm}\label{per-thm1}
Let $0<i_k<n_k$, and let $\omega_{i_0\ldots i_d}$ 
be the differential forms in \S \ref{form-sect}.
Put $a_k:=1-i_k/n_k$.
Then there is a homology cycle $\Delta_\alpha\in H_d(U_\alpha,\Z)$ which is a
vanishing cycle at $\alpha=0$ such that
\begin{equation}\label{per-thm1-eq0}
\int_{\Delta_\alpha} \omega_{i_0\cdots i_d}=\frac{(2\pi i)^d}{n_0\cdots n_d}\cdot
{}_{d+1}F_d\left({a_0,\ldots,a_d\atop1,\ldots,1};\alpha\right).
\end{equation}
\end{thm}
\begin{pf}
We define the cycle $\Delta_\alpha$ to be the following.
Recall the equation $(1-x_0^{n_0})\cdots(1-x_d^{n_d})=\alpha$ of $U_\alpha$.
Let $0<|\alpha|\ll1$.
When $|x_1-1|=\cdots=|x_d-1|=|\alpha|^{\frac{1}{d+1}}$, then $x_0$ such that $|x_0-1|\ll1$
is uniquely determined.
Hence a torus ${\mathbb T}=\{(x_1,\ldots,x_d)\mid |x_1-1|=\cdots=|x_d-1|=|\alpha|^{\frac{1}{d+1}}\}$
defines a homology cycle $\Delta_\alpha$ in $U_\alpha$
with a suitable orientation.
By definition, $\Delta_\alpha$ is a vanishing cycle at $\alpha=0$.
We show \eqref{per-thm1-eq0}.
\begin{align*}
\int_{\Delta_\alpha} \omega_{i_0\ldots i_d}
&=n_0^{-1}\int_{\Delta_\alpha}x_0^{i_0-n_0}x_1^{i_1-1}\cdots x_d^{i_d-1}
\frac{dx_1\cdots dx_d}{(1-x_1^{n_1})\cdots(1-x_d^{n_d})}\\
&=n_0^{-1}\int_{\Delta_\alpha}
\left(1-\frac{\alpha}{(1-x_1^{n_1})\cdots(1-x_d^{n_d})}\right)^{-a_0}
\frac{x_1^{i_1-1}\cdots x_d^{i_d-1}dx_1\cdots dx_d}{(1-x_1^{n_1})\cdots(1-x_d^{n_d})}
\\
&=n_0^{-1}\int_{\Delta_\alpha}\left(
\sum_{n=0}^\infty
\frac{(a_0)_n}{n!}\frac{x_1^{i_1-1}\cdots x_d^{i_d-1}dx_1\cdots dx_d}{(1-x_1^{n_1})^{n+1}
\cdots(1-x_d^{n_d})^{n+1}}\alpha^n\right)
\\
&=n_0^{-1}
\sum_{n=0}^\infty\frac{(a_0)_n}{n!}\alpha^n
\int_{\Delta_\alpha}
\frac{x_1^{i_1-1}\cdots x_d^{i_d-1}dx_1\cdots dx_d}{(1-x_1^{n_1})^{n+1}
\cdots(1-x_d^{n_d})^{n+1}}
\\
&=n_0^{-1}
\sum_{n=0}^\infty\frac{(a_0)_n}{n!}\alpha^n
\prod_{k=1}^d\oint_{|x_k-1|=\ve}
\frac{x_k^{i_k-1}dx_k}{(1-x_k^{n_k})^{n+1}}
\end{align*}
where the interchange of the integral and summation can be verified due to 
the uniform convergence by the assumption $|\alpha|\ll1$.
Applying Lemma \ref{lem-per-thm1-eq2}, \eqref{per-thm1-eq0} follows.
\end{pf}
\begin{cor}\label{per-cor1}
Let $\wt\omega^{(r)}_{i_0\cdots i_d}\in H^d_\dR(X/A)$ be the lifting 
in Lemma \ref{form-lem2}. Then
\[
\int_{\Delta_\alpha} \wt\omega^{(r)}_{i_0\cdots i_d}=\frac{(2\pi i)^d}{n_0\cdots n_d}
(a_0(a_0-1)\cdots(a_0-r+1))^{-1}\alpha^r F^{(r)}(\alpha)
\]
where
\[
F^{(r)}(t):=\frac{d^r}{dt^r}\left({}_{d+1}F_d\left({a_0,\ldots,a_d\atop1,\ldots,1};t\right)\right).
\]
\end{cor}
\begin{pf}
Immediate from Theorem \ref{per-thm1} and Lemma \ref{form-lem3}.
\end{pf}
\subsection{Cohomology of Hypergeometric Schemes}\label{dim-sect}
Let $\mu_m=\mu_m(\C)$ denote the group of $m$-th roots of unity in $\C$.
A finite group $G=\mu_{n_i}\times \cdots\times \mu_{n_d}$ acts on $U$ and $U_\alpha$
as in \S \ref{setting-sect}.
The topological fundamental group
$\pi_1(S,\alpha)$ acts on the Betti cohomology groups
$H^\bullet_B(X_\alpha,k)$ and $H^\bullet_B(U_\alpha,k)$ where $k$ is a commutative ring,
and it commutes with the actionof $G$.
The cohomology groups $H^\bullet_B(U_\alpha,\Q)$ carry the mixed Hodge structures
by Deligne. We denote by
 $W_* H^\bullet_B(U_\alpha,\Q)$ and $W_* H^\bullet_\dR(U_\alpha/\C)$ 
the weight filtrations.

\medskip

Let $a_0,\ldots,a_d\in \C$.
Let
\begin{equation}\label{per-lem1-eq1}
\phi:\pi_1(S,\alpha)\lra \mathrm{GL} (V_\HG(a_0,\ldots,a_d)_\alpha)
\end{equation}
be the monodromy representation of the hypergeometric function
\[
{}_{d+1}F_d\left({a_0,\ldots,a_d\atop1,\ldots,1};t\right).
\]
This is defined in the following way.
Let $V_\HG(a_0,\ldots,a_d)_\alpha$ be the complex linear subspace in the stalk
$\O^{\text{an}}_{S,\alpha}\cong\C\{t-\alpha\}$ which is generated by
all analytic continuations of the above function.
The action of $\pi_1(S,\alpha)$ on the space $V_\HG(a_0,\ldots,a_d)_\alpha$
is defined in a natural way.

The linear space $V_\HG(a_0,\ldots,a_d)_\alpha$ is at most 
$(d+1)$-dimensional.
We denote by $\cV_\HG(a_0,\ldots,a_d)$ the corresponding locally constant sheaf.
It is a fundamental theorem that 
if $a_i\not\in\Z$ for all $i$, then 
\eqref{per-lem1-eq1} is a $(d+1)$-dimensional irreducible
representation (\cite[Proposition 3.3]{BH}).

\begin{lem}\label{per-lem1}
Suppose that $n_k>1$ for all $k$.
Let $(i_0,\ldots, i_d)$ be a $(d+1)$-tuple of inetgers such that $0<i_k<n_k$. 
We denote by $V(i_0,\ldots,i_d)$ the eigenspace as in \S \ref{setting-sect}.
Put $a_k:=1-i_k/n_k$.
Let $V_{\dR,\alpha}\subset W_dH^d_\dR(U_\alpha/\C)(i_0,\ldots,i_d)$ be the subspace generated by
$\wt\omega^{(0)}_{i_0\ldots i_d},\ldots,\wt\omega^{(d)}_{i_0\ldots i_d}$ 
in Lemma \ref{form-lem3}.
Then \[
\dim_\C V_{\dR,\alpha}=d+1.
\]
If $\dim_\C W_dH^d_\dR(U_\alpha/\C)(i_0,\ldots,i_d)=d+1$
namely $V_{\dR,\alpha}= W_dH^d_\dR(U_\alpha/\C)(i_0,\ldots,i_d)$, then 
\[
W_dH^d_B(U_\alpha,\C)(i_0,\ldots,i_d)\cong V_\HG(a_0,\ldots,a_d)^\vee_\alpha
\]
as $\C[\pi_1(S,\alpha)]$-module where $(-)^\vee$ denotes the dual representation. 
In particular, this is irreducible.
\end{lem}
We shall soon see 
that $\dim W_dH^d_\dR(U_\alpha/\C)(i_0,\ldots,i_d)=d+1$ always holds (Theorem \ref{dim-thm1}).
\begin{pf}
We write $V_{\HG}=V_\HG(a_0,\ldots,a_d)$ simply.
Let $\Delta_\alpha\in H_d(U_\alpha,\Z)$ be the homology cycle in Theorem \ref{per-thm1}.
Let $V_{B,\alpha}\subset H_d(U_\alpha,\C)$ be the sub $\C[\pi_1(S,\alpha)]$-module generated
by $\Delta_\alpha$.
There is a $\C$-linear map
\[
\Phi:V_{B,\alpha}\lra V_{\mathrm{HG},\alpha}
,\quad \gamma\longmapsto
\int_\gamma\omega_{i_0\ldots i_d}
\]
which is compatible with respect to the action of $\pi_1(S,\alpha)$.
Since $V_{\text{HG},\alpha}$ is irreducible,
$\Phi$ is surjective.
This implies that there are homology cycles $\gamma_0,\ldots,\gamma_d\in V_{B,\alpha}$ 
such that $\int_{\gamma_j}\omega_{i_0\ldots i_d}=F_j(t)$.
By Corollary \ref{per-cor1}, we have
\[
\int_{\gamma_j}\wt\omega^{(r)}_{i_0\ldots i_d}=2\pi \sqrt{-1}(n_0\cdots n_d)^{-1}t^rF_j^{(r)}(t).
\]
Now the linear independence of $\wt\omega^{(0)}_{i_0\ldots i_d},\ldots,
\wt\omega^{(d)}_{i_0\ldots i_d}$ follows from the non-vanishing of the Wronskian
determinant
\[
\det\begin{pmatrix}
F_0(\alpha)&\cdots&F_d(\alpha)\\
F_0^{(1)}(\alpha)&\cdots&F^{(1)}_d(\alpha)\\
\vdots&&\vdots\\
F_0^{(d)}(\alpha)&\cdots&F^{(d)}_d(\alpha)
\end{pmatrix}\ne0,\quad \forall\,\alpha\in\C\setminus\{0,1\}.
\]
We show the latter assertion.
The map $\Phi$ factors through the eigenspace $V_{B,\alpha}(i_0,\ldots,i_d)$,
and also the image $\ol{V}_{B,\alpha}(i_0,\ldots,i_d)
:=\Image[V_{B,\alpha}(i_0,\ldots,i_d)\to H_d(X_\alpha,\C)(i_0,\ldots,i_d)]$
as $\omega_{i_0\ldots i_d}\in \vg(X_\alpha,\Omega^d_{X_\alpha})$,
\[
\xymatrix{
V_{B,\alpha}(i_0,\ldots,i_d)\ar[r]^{\text{surj.}}&
\ol{V}_{B,\alpha}(i_0,\ldots,i_d)\ar[r]^{\text{surj.}}\ar[d]_\bigcap&
V_{\mathrm{HG},\alpha}\\
&\Hom(W_dH^d_\dR(U_\alpha/\C)(i_0,\ldots,i_d),\C).
}
\]
Suppose that $\dim_\C W_dH^d_\dR(U_\alpha/\C)(i_0,\ldots,i_d)=d+1$.
Since $V_\HG$ is an irreducible $(d+1)$-dimensional representation,
the above diagram implies isomorphisms
\[
(W_dH^d_B(U_\alpha,\C)(i_0,\ldots,i_d))^\vee\cong
\ol{V}_{B,\alpha}(i_0,\ldots,i_d)\cong
V_{\mathrm{HG},\alpha}
\]
of $\pi_1(S,\alpha)$-modules.
This completes the proof.
\end{pf}

\begin{thm}\label{dim-thm1}
Let $n_k\geq 1$ be integers.
\begin{enumerate}
\renewcommand{\labelenumi}{\rm{(\roman{enumi})}}
\item
Let $U_1^\sm=\Spec \C[x_0,\ldots,x_d]/((1-x_0^{n_0})\cdots(1-x_d^{n_d})-1)
\setminus\{(0,\ldots,0)\}$.
Then $W_iH^i(U_1^\sm,\Q)=0$ if $1\leq i\leq d-1$.
\item
Let $U_\alpha=\Spec \C[x_0,\ldots,x_d]/((1-x_0^{n_0})\cdots(1-x_d^{n_d})-\alpha)$
with $\alpha\ne0,1$.
Then $W_iH^i(U_\alpha,\Q)=0$ if $1\leq i\leq d-1$.
\item
If $n_k=1$ for some $k$, then $W_dH^d(U_\alpha,\C)=0$.
Suppose that $n_k>1$ for all $k$.
Then
\[
W_dH^d(U_\alpha,\C)=\bigoplus_{i_0,\ldots,i_d}W_dH^d(U_\alpha,\C)(i_0,\ldots,i_d)
\]
and $\dim W_dH^d(U_\alpha,\C)(i_0,\ldots,i_d)=d+1$ where
$(i_0,\ldots,i_d)$ runs over all $(d+1)$-tuple of integers such that $0<i_k<n_k$.
\end{enumerate}
\end{thm}
\begin{pf}
In case that $n_k=1$ for some $k$, one can easily prove (i), (ii) and (iii)
on noticing that
\begin{align*}
U_\alpha
&\cong
\Spec \C[x_0,\ldots,\wh{x_k},\ldots,x_d][(1-x_0^{n_0})^{-1}\cdots\wh{(1-x_k)}^{-1}\cdots(1-x_d^{n_d})^{-1}]\\
&\cong\prod_{0\leq i\leq d,\,i\ne k}\bA^1(x_i)\setminus\{x_i^{n_i}=1\}
\end{align*}
where $\alpha\in \C\setminus\{0\}$ (including $\alpha=1$).

Suppose that $n_k\geq2$ for all $k$.
We show (i), (ii) and (iii) by the induction on $d$.
We denote by (i)$_d$, (ii)$_d$ and (iii)$_d$ the statements for $d$.
There is nothing to prove for (i)$_1$ and (ii)$_1$.
We show (iii)$_1$.
Let $X_\alpha\supset U_\alpha$ be the smooth compactification.
It follows from Lemma \ref{per-lem1} that
$\dim H^1(X_\alpha,\C)(i_0,i_1)=\dim W_1H^1(U_\alpha,\C)(i_0,i_1)\geq 2$. 
Since the genus of $X_\alpha$ is $(n_0-1)(n_1-1)$, we have
\[
2(n_0-1)(n_1-1)\leq \sum_{i_0,i_1}\dim H^1(X_\alpha,\C)(i_0,i_1)=\dim H^1(X_\alpha,\C)=2(n_0-1)(n_1-1)
\]
and hence the equality holds, which implies
$\dim W_1H^1(U_\alpha,\C)(i_0,i_1)=2$ for all $i_0,i_1$ as required. 

\medskip

Let $d\geq 2$. Suppose that (i)$_{d-1}$, (ii)$_{d-1}$ and (iii)$_{d-1}$ hold.
We first show (i)$_{d}$.
Let $\bar S:=\bA^1(x_d)\setminus\{x_d^{n_d}=1\}\supset S:=\bA^1(x_d)\setminus\{x_d^{n_d}=1,0\}$.
Let $g:U^\sm_1\to \bar S$ be the projection given by $(x_0,\ldots,x_d)\mapsto x_d$,
and put $U_1^{\sm,\circ}:=g^{-1}(S)$,
\begin{equation}\label{dim-thm1-eq1}
\xymatrix{
U_1^\sm\ar[r]^g& \bar S\\
U_1^{\sm,\circ}\ar[r]\ar[u]&S.\ar[u]
}
\end{equation}
Since $U_1^{\sm,\circ}\to S$ is a topological fibartion, the sheaves $R^ig_*\Q|_S$ are
locally constant sheaves.
Therefore one has
\begin{equation}\label{dim-thm1-eq2}
H^i(U_1^{\sm,\circ},\Q)=H^0(S,R^ig_*\Q)\op H^1(S,R^{i-1}g_*\Q),
\end{equation}
and hence
\[
W_iH^i(U_1^{\sm,\circ},\Q)=H^0(S,W_iR^ig_*\Q)\op W_iH^1(S,R^{i-1}g_*\Q)
\]
for all $i\geq 0$. The map
 $W_iH^1(S,W_{i-1}R^{i-1}g_*\Q)\to W_iH^1(S,R^{i-1}g_*\Q)$
is surjective. By (ii)$_{d-1}$, one has the vanishing $W_jR^jg_*\Q=0$ for all $0<j<d-1$.
Therefore, 
$W_iH^i(U_1^{\sm,\circ},\Q)=0$ if $2\leq i\leq d-2$.
If $i=1$, one also has $W_1H^1(U_1^{\sm,\circ},\Q)=W_iH^1(S,\Q)=0$.
Let $i=d-1$. 
One has
\begin{equation}\label{dim-thm1-eq2-d}
W_{d-1}H^{d-1}(U_1^{\sm,\circ},\Q)=H^0(S,W_{d-1}R^{d-1}g_*\Q).
\end{equation}
Let $\cV_\HG(a_0,\ldots,a_{d-1})$ be the locally constant sheaf corresponding 
to the monodromy representation \eqref{per-lem1-eq1} of the
hypergeometric functions
\[
{}_dF_{d-1}\left({a_0,\ldots,a_{d-1}\atop1,\ldots,1};t\right),\quad a_k:=1-\frac{i_k}{n_k}.
\]
By (iii)$_{d-1}$ and Lemma \ref{per-lem1}, 
$W_{d-1}R^{d-1}g_*\C$ is the direct sum of 
$W_{d-1}R^{d-1}g_*\C(i_0,\ldots,i_{d-1})$'s, and each $W_{d-1}R^{d-1}g_*\C(i_0,\ldots,i_{d-1})$
is isomorphic to the dual of $\rho^*\cV_\HG(a_0,\ldots,a_{d-1})$
where $\rho:S\to\bA^1(t)\setminus\{0,1\}$ is the morphism such that 
$\rho^*(t)=1/(1-x_d^{n_d})$.
In particular, the right hand side of \eqref{dim-thm1-eq2-d}
vanishes. We now have 
\begin{equation}\label{dim-thm1-ind-eq1}
W_iH^i(U_1^{\sm,\circ},\Q)=0,\quad 1\leq\forall\, i\leq d-1.
\end{equation}
We show the vanishing $W_iH^i(U^\sm_1,\Q)=0$ for $1\leq i\leq d-1$.
Put $D^\sm_1:=U_1^\sm\setminus U^{\sm,\circ}_1=
\Spec\C[x_0,\ldots,x_{d-1}]/
((1-x_0^{n_0})\cdots(1-x_{d-1}^{n_{d-1}})-1)\setminus\{(0,\ldots,0)\}$
a smooth connected divisor.
There is the localization exact sequence
\begin{equation}\label{dim-thm1-eq3}
\cdots\to
H^{i-1}(U^{\sm,\circ}_1,\Q)\to 
H^{i-2}(D^\sm_1,\Q(-1))\to 
H^{i}(U_1^\sm,\Q)\to H^{i}(U^{\sm,\circ}_1,\Q) \to\cdots
\end{equation}
of mixed Hodge structures, which gives rise to an exact sequence
\[
W_{i-2}H^{i-2}(D^\sm_1,\Q)\lra
W_iH^{i}(U_1^\sm,\Q)\lra W_iH^{i}(U^{\sm,\circ}_1,\Q).
\]
By \eqref{dim-thm1-ind-eq1} the right term vanishes if $1\leq i\leq d-1$.
The left term vanishes if $1\leq i\leq d$ and $i\ne2$ by (i)$_{d-1}$.
In case $i=2$, since $H^1(U^{\sm,\circ}_1,\Q)\to H^0(D_1,\Q(-1))=\Q$ is surjective,
we also have the vanishing of the middle term.
This competes the proof of (i)$_d$.

\medskip

Next we show (ii)$_d$.
Let $g_\alpha:U_\alpha\to \bar S$ be the projection given by $(x_0,\ldots,x_d)\mapsto x_d$,
and put $S_\alpha:=\bA^1(x_d)\setminus\{x_d^{n_d}=1,1-\alpha\}$
and $U_\alpha^\circ:=g^{-1}_\alpha(S)$,
\begin{equation}\label{dim-thm1-eq4}
\xymatrix{
U_\alpha\ar[r]^{g_\alpha}& \bar S\\
U_\alpha^\circ\ar[r]\ar[u]&S_\alpha.\ar[u]
}
\end{equation}
We note that $U_\alpha^\circ\to S_\alpha$ is a topological fibration.
In the same way as the proof of \eqref{dim-thm1-ind-eq1}, one can show the vanishing
\begin{equation}\label{dim-thm1-ind-eq2}
W_iH^i(U_\alpha^\circ,\Q)=0,\quad 1\leq i\leq d-1.
\end{equation}
For an integer $1\leq i\leq n_d$, let $D_i$ be the fiber of $g_\alpha$ at $x_d=\sqrt[n_d]{1-\alpha}\zeta_{n_d}^i$
where $\zeta_{n_d}$ is a fixed primitive $n_d$-th root of unity.
The divisor $D_i$ has a unique singular point $z_i=(0,\ldots,0,
\sqrt[n_d]{1-\alpha}\zeta_{n_d}^i)$.
Put $D_i^\sm:=D_i\setminus\{z_i\}$ a smooth connected divisor.
Put $D:=\coprod_{i=1}^{n_d} D_i$, $D^\sm:=\coprod_{i=1}^{n_d} D^\sm_i$ 
and $Z:=\{z_1,\ldots,z_{n_d}\}$.
There is the localization exact sequence
\begin{equation}\label{dim-thm1-eq5}
\cdots\to
H^{i-1}(U^{\circ}_\alpha,\Q)\to 
H^{i-2}(D^\sm,\Q(-1))\to 
H^{i}(U_\alpha\setminus Z,\Q)\to H^{i}(U^{\circ}_\alpha,\Q) \to\cdots
\end{equation}
of mixed Hodge structures,
and this gives rise to an exact sequence
\begin{equation}\label{dim-thm1-ind-eq3}
W_{i-2}H^{i-2}(D^\sm,\Q)\lra
W_iH^{i}(U_\alpha\setminus Z,\Q)\lra W_iH^{i}(U^{\circ}_\alpha,\Q).
\end{equation}
By \eqref{dim-thm1-ind-eq2}
the right term vanishes if 
$1\leq i\leq d-1$.
The left term vanishes if $1\leq i\leq d$ and $i\ne2$ by (i)$_{d-1}$.
In case $i=2$, we also have the vanishing of the middle term
as
$H^1(U^\circ_\alpha,\Q)\to H^0(D^\sm,\Q(-1))=\Q^{\op n_d}$ is surjective.
We thus have
\[
W_iH^i(U_\alpha\setminus Z,\Q)=0,\quad 1\leq \forall\,i\leq d-1.
\]
One can replace $U_\alpha\setminus Z$ with $U_\alpha$ in the above
as $H^i(U_\alpha,\Q)\cong
H^i(U_\alpha\setminus Z,\Q)$ for $i\ne 2d,2d-1$.
This completes the proof of (ii)$_d$.

\medskip

Finally we show (iii)$_d$.
By Lemma \ref{per-lem1}, the inequality
$\dim W_dH^d(U_\alpha,\C)(i_0,\ldots,i_d)\geq d+1$ always holds.
Therefore it is enough to show 
\begin{equation}\label{dim-thm1-ind-eq4}
\dim W_dH^d(U_\alpha,\Q)\leq (d+1)(n_0-1)\cdots (n_d-1).
\end{equation}
Let $i=d$ in the exact sequence \eqref{dim-thm1-ind-eq3}.
If $d\geq 3$, the left term vanishes by (i)$_{d-1}$.
If $d=2$, then $H^1(U^\circ_\alpha,\Q)\to H^0(D^\sm,\Q(-1))=\Q^{\op n_d}$ is surjective.
In both cases, the map
\[
W_dH^d(U_\alpha,\Q)=
W_dH^d(U_\alpha\setminus Z,\Q)\lra W_dH^d(U^{\circ}_\alpha,\Q)
\]
is injective.
Noticing that $U_\alpha^\circ\to S_\alpha$ is a topological fibration and fibers
are smooth affine varieties of dimension $d-1$,
one has
\begin{equation}\label{dim-thm1-eq6}
H^d(U^{\circ}_\alpha,\Q)=H^0(S_\alpha, R^dg_{\alpha*}\Q)\op 
H^1(S_\alpha, R^{d-1}g_{\alpha*}\Q)
=H^1(S_\alpha, R^{d-1}g_{\alpha*}\Q)
\end{equation}
where the vanishing $R^dg_{\alpha*}\Q=0$ follows by the Lefschetz affine theorem.
We have
\[
\dim W_dH^d(U_\alpha, \Q)\leq
\dim W_dH^1(S_\alpha, R^{d-1}g_{\alpha*}\Q)
\leq\dim W_dH^1(S_\alpha, W_{d-1}R^{d-1}g_{\alpha*}\Q).
\]
We show the last term is $(d+1)(n_0-1)\cdots (n_d-1)$, which implies 
\eqref{dim-thm1-ind-eq4}.
Put $\cV=W_{d-1}R^{d-1}g_{\alpha*}\C$.
By (iii)$_{d-1}$, one has
\[
\cV=\bigoplus_{i_0,\ldots,i_{d-1}}\cV(i_0,\ldots,i_{d-1}),\quad \dim\cV(i_0,\ldots,i_{d-1})=d
\]
where $(i_0,\ldots,i_{d-1})$ runs over all $d$-tuple of integers such that $0<i_k<n_k$.
By Lemma \ref{per-lem1},
each $\cV(i_0,\ldots,i_{d-1})$ is isomorphic to the dual of 
$\rho_\alpha^*\cV_\HG(a_0,\ldots,a_{d-1})$
where $\rho_\alpha:S_\alpha\to\bA^1(t)\setminus\{0,1\}$ is the morphism such that
$\rho^*(t)=\alpha/(1-x_d^{n_d})$.
In particular $\vg(S_\alpha,\cV)=0$.
We have
\[
\dim H^1(S_\alpha,\cV)=-\chi(S_\alpha,\cV)
=-\chi(S_\alpha)\dim\cV=(2n_d-1)\cdot d(n_0-1)\cdots (n_{d-1}-1).
\]
Let $j:S_\alpha\hra\P^1$ be the open immersion.
There is an exact sequence
\[
\xymatrix{
0\ar[r]&
H^1(\P^1,j_*\cV)\ar[r]&
H^1(S_\alpha,\cV)\ar[r]&
\vg(\P^1,R^1j_*\cV)\ar[r]&0\\
&W_dH^1(S_\alpha,\cV)\ar@{=}[u]
}\]
where the equality follows by the fact that the weight of $R^1j_*\cV$ is $\geq d+1$.
The sheaf $R^1j_*\cV$ is supported on 
at most $x_d^{n_d}=1,1-\alpha,\infty$.
Let $T_p\in\pi_1(S,\alpha)$ be the local monodromy at $p$.
Then
\[
(R^1j_*\cV(i_0,\ldots,i_{d-1}))_{x_d=\zeta_{n_d}^i}=\Coker[T_\infty-\id:
V_\HG(a_0,\ldots,a_{d-1})_\alpha\to V_\HG(a_0,\ldots,a_{d-1})_\alpha]=0.
\]
Moreover we have that
\[
(R^1j_*\cV(i_0,\ldots,i_{d-1}))_{x_d=\sqrt[n_d]{1-\alpha}\zeta_{n_d}^i}
=\Coker[T_1-\id:
V_\HG(a_0,\ldots,a_{d-1})_\alpha\to V_\HG(a_0,\ldots,a_{d-1})_\alpha]
\]
is $(d-1)$-dimensional, and 
\[
(R^1j_*\cV(i_0,\ldots,i_{d-1}))_{x_d=\infty}
=\Coker[T_0^{n_d}-\id:
V_\HG(a_0,\ldots,a_{d-1})_\alpha\to V_\HG(a_0,\ldots,a_{d-1})_\alpha
\]
is $1$-dimensional
(note that $T_0$ is unipotent).
We thus have 
\begin{align*}
\dim W_dH^1(S_\alpha,\cV)
&=
d(2n_d-1)\cdot (n_0-1)\cdots (n_{d-1}-1)\\
&\hspace{2cm}-
((d-1)n_d+1)\cdot (n_0-1)\cdots (n_{d-1}-1)\\
&=
(d+1) (n_0-1)\cdots (n_d-1)
\end{align*}
as required.
This completes the proof of \eqref{dim-thm1-ind-eq4}, and hence (iii)$_d$.
\end{pf}
\begin{cor}\label{dim-cor}
Suppose $n_k>1$ for all $k$. Let $0<i_k<n_k$ and put $a_k:=1-i_k/n_k$.
Then the representation $W_dH^d(U_\alpha,\C)(i_0,\ldots,i_d)$ of $\pi_1(\bA^1(t)\setminus\{0,1\},\alpha)$
is isomorphic to the dual of the
monodromy representation of 
\[
{}_{d+1}F_d\left({a_0,\ldots,a_d\atop1,\ldots,1};t\right).
\]
\end{cor}
\begin{pf}
This is immediate from Lemma \ref{per-lem1} together with 
Theorem \ref{dim-thm1} (iii).
\end{pf}
\begin{cor}\label{dim-cor2}
Let $\cD:=\C\langle t,(t-t^2)^{-1},\frac{d}{dt}\rangle$ be the Wyle algebra, which 
acts on $W_dH^d_\dR(U/A)$.
Put $D:=t\frac{d}{dt}$ and $P_\HG:=D^{d+1}-t(D+a_0)\cdots(D+a_d)$.
Let $\omega_{i_0\ldots i_d}\in \vg(X,\Omega^d_{X/A})$ be the regular $d$-form
\eqref{form-eq1}. 
Then there is an isomorphism
\[
\cD/\cD P_\HG\os{\cong}{\lra}
W_dH^d_\dR(U/A)(i_0,\ldots,i_d),\quad P\longmapsto P(\omega_{i_0\ldots i_d})
\]
of $\cD$-modules.
\end{cor}
\begin{pf}
By Theorem \ref{dim-thm1} (iii) and Lemma \ref{per-lem1},
$W_dH^d_\dR(U/A)(i_0,\ldots,i_d)$ is a free $A$-module with
basis $\wt\omega^{(0)}_{i_0\ldots i_d},\ldots,\wt\omega^{(d)}_{i_0\ldots i_d}$.
Let $F_0(t),\ldots,F_d(t)\in \O_{S,\alpha}$ be analytic continuations of
${}_{d+1}F_d\left({a_0,\ldots,a_d\atop1,\ldots,1};t\right)$
which are linearly independent over $\C$. They are the solutions of
$P_\HG$.
Then, as is shown in the proof of Lemma \ref{per-lem1},
the map 
\[
W_dH^d(U_\alpha,\C)(i_0,\ldots,i_d)\lra (\O_{S,\alpha})^{\op d+1},\quad
\wt\omega^{(r)}_{i_0\ldots i_d}\longmapsto (F_0^{(r)}(t),\ldots,F_d^{(r)}(t))
\]
is injective, so that we have
\[
P(\omega_{i_0\ldots i_d})=0\quad \Longleftrightarrow\quad
P(F_i(t))=0,\,\forall\,i
\]
for $P\in \cD$. In particular a vanishing
$P_\HG(\omega_{i_0\ldots i_d})=0$ follows.
We thus have a homomorphism
\[
\cD/\cD P_\HG\lra
W_dH^d_\dR(U/A)(i_0,\ldots,i_d),\quad P\longmapsto P(\omega_{i_0\ldots i_d})
\]
of $\cD$-modules.
Since both are irreducible smooth $\cD$-module with regular singularities, this is bijective (Schur's lemma). 
\end{pf}
\begin{thm}\label{dim-thm2}
Suppose $n_k>1$ for all $k$. 
Let $(i_0,\ldots,i_d)$ be $(d+1)$-tuple of integers such that $0< i_k<n_k$ for all $k$.
\begin{enumerate}
\renewcommand{\labelenumi}{\rm{(\roman{enumi})}}
\item
$H^i(U_1^\sm,\C)(i_0,\ldots,i_d)=0$ for all $i<d$.
\item
$H^i(U_\alpha,\C)(i_0,\ldots,i_d)=0$ for all $i<d$.
\item
$H^d(U_\alpha,\C)(i_0,\ldots,i_d)=W_dH^d(U_\alpha,\C)(i_0,\ldots,i_d) $.
\end{enumerate}
\end{thm}
\begin{pf}
We show the theorem by the induction on $d$.
We denote by (i)$_d$, (ii)$_d$ and (iii)$_d$ the statements for $d$.
We show the case $d=1$.
It is simple to see (i)$_1$ and (ii)$_1$.
We show that (iii)$_1$ holds.
Let $X_\alpha\supset U_\alpha$ be the smooth compactification, and put $Z:=
X_\alpha\setminus U_\alpha$.
Then there is the exact sequence
\[
0\lra H^1(X_\alpha,\Q)\lra H^1(U_\alpha,\Q)\lra H^0(Z,\Q)
\]
which is compatible with the actions of $\sigma_i(\nu_i)$'s.
Since either of $\sigma_0(\nu_0)$ or $\sigma_1(\nu_1)$ acts on a point $z\in Z$
as identity, this implies $H^1(U_\alpha,\C)(i_0,i_1)\subset H^1(X_\alpha,\C)(i_0,i_1)$
for any $0<i_k<n_k$. This completes the proof of (iii)$_1$.

\medskip

Let $d\geq 2$. Suppose that (i)$_{d-1}$, (ii)$_{d-1}$ and (iii)$_{d-1}$ hold.
We first show (i)$_d$.
We consider the diagram \eqref{dim-thm1-eq1},
\[
\xymatrix{
U_1^\sm\ar[r]^{g\hspace{1.7cm}}& \bar S=\bA^1(x_d)\setminus\{x_d^{n_d}=1\}\\
U_1^{\sm,\circ}\ar[r]\ar[u]&S=\bA^1(x_d)\setminus\{x_d^{n_d}=0,1\}.\ar[u]
}
\]
Put $D^\sm_1:=U_1^\sm\setminus U^{\sm,\circ}_1$ a smooth connected divisor.
The localization exact sequence \eqref{dim-thm1-eq3}
is compatible with the actions of the automorphisms
$\sigma_i(\nu_i)$ in \eqref{setting-eq1}.
Since $\sigma_d(\nu_d)$ acts on $D^\sm_1$ as the identity, 
one has 
\begin{equation}\label{dim-thm2-eq1}
H^i(U_1^\sm,\C)(i_0,\ldots,i_d)\os{\cong}{\lra} H^{i}(U^{\sm,\circ}_1,\C)(i_0,\ldots,i_d),
\quad \forall\,i\geq 0.
\end{equation}
The isomorphism \eqref{dim-thm1-eq2} yields
\[
H^i(U_1^{\sm,\circ},\C)(i_0,\ldots,i_d)
\subset H^0(S,R^ig_*\C(i_0,\ldots,i_{d-1}))\op H^1(S,R^{i-1}g_*\C(i_0,\ldots,i_{d-1}))
\]
for all $i$.
By (ii)$_{d-1}$ one has $R^jg_*\C(i_0,\ldots,i_{d-1})=0$ for $j<d-1$.
Therefore the right hand side vanishes if $i<d-1$, and hence one has
$H^i(U_1^\sm,\C)(i_0,\ldots,i_d)=0$ for all $i<d-1$ by \eqref{dim-thm2-eq1}.
Let $i=d-1$. Then
\[
H^{d-1}(U_1^{\sm},\C)(i_0,\ldots,i_d)=
H^{d-1}(U_1^{\sm,\circ},\C)(i_0,\ldots,i_d)
\subset H^0(S,R^{d-1}g_*\C(i_0,\ldots,i_{d-1})).
\]
By (iii)$_{d-1}$ together with Corollary \ref{dim-cor}, the most right term vanishes.
This completes the proof of (i)$_d$.

\medskip

We show (ii)$_d$.
Let 
\[
\xymatrix{
U_\alpha\ar[r]^{g_\alpha\hspace{1.7cm}}& \bar S=\bA^1(x_d)\setminus\{x_d^{n_d}=1\}\\
U_\alpha^\circ\ar[r]\ar[u]&S_\alpha=\bA^1(x_d)\setminus\{x_d^{n_d}=1,1-\alpha\}.
\ar[u]
}
\]
be the diagram \eqref{dim-thm1-eq4}. In the same way as above one can show
the vanishing 
\[
H^i(U_\alpha^\circ,\C)(i_0,\ldots,i_d)=0,\quad i<d.
\]
We use the notation $Z\subset D=U_\alpha\setminus U_\alpha^\circ$ and
$D^\sm=D\setminus Z$ after \eqref{dim-thm1-ind-eq2}.
The localization sequence \eqref{dim-thm1-eq5}
is compatible with the actions of $\sigma_i(\nu_i)$'s.
By (i)$_{d-1}$, one has $H^j(D^\sm,\C)(i_0,\ldots,i_{d-1})=0$ for all $j<d-1$.
Hence
\[
H^i(U_\alpha\setminus Z,\C)(i_0,\ldots,i_d)\os{\cong}{\lra}
H^i(U_\alpha^\circ,\C)(i_0,\ldots,i_d)=0,\quad i<d.
\]
Now the assertion (ii)$_d$ follows from the fact $H^i(U_\alpha,\Q)
=H^i(U_\alpha\setminus Z,\Q)$ for $i\ne 2d,2d-1$.

\medskip

Finally we show (iii)$_d$.
It is enough to show that 
$H^d(U_\alpha,\C)(i_0,\ldots,i_d)$ has pure weight $d$.
By an exact sequence
\[
\xymatrix{
H^{d-2}(D^\sm)(i_0,\ldots,i_d)\ar[r]&
H^d(U_\alpha\setminus Z)(i_0,\ldots,i_d)\ar[r]&
H^d(U_\alpha^\circ,\C)(i_0,\ldots,i_d)\\
0\ar@{=}[u]_{\text{by (i)}_{d-1}}&H^d(U_\alpha)(i_0,\ldots,i_d)\ar[u]_\cong
}\]
together with the isomorphism \eqref{dim-thm1-eq6},
one has
\[
H^d(U_\alpha,\C)(i_0,\ldots,i_d)\subset
H^1(S_\alpha,\cM(i_0,\ldots,i_{d-1})).
\]
where we put $\cM:=R^{d-1}g_{\alpha*}\C$.
Let $j:S_\alpha\hra \P^1(x_d)$ be the open immersion.
Then there is a commutative diagram
\[
\xymatrix{
0\ar[r]&H^1(\P^1,j_*\cM)
\ar[r]&
H^1(S_\alpha,\cM)\ar[r]&
H^0(\P^1,R^1j_*\cM)\\
&&H^d(U_\alpha,\C)\ar[r]^{h\hspace{1.7cm}}
\ar[u]^\cup&\bigoplus_{\beta=1,1-\alpha,\infty} (R^1j_*\cM)_{x_d^{n_d}=\beta}\ar@{=}[u]
}\]
with exact row.
The cohomology $H^1(\P^1,j_*\cM)$ has pure weight $d$.
Therefore it is enough to show that the image $h(H^d(U_\alpha,\C)(i_0,\ldots,i_d))$ vanishes.
The image of $h$ lies at most in the components of $\beta=1,\infty$.
Since $\sigma_d(\nu_d)$ acts on the component
of $x_d=\infty$ as identity, the image of $H^d(U_\alpha,\C)(i_0,\ldots,i_d)$
lies only in the component of $x_d^{n_d}=1$,
\[
h(H^d(U_\alpha,\C)(i_0,\ldots,i_d))\subset
\bigoplus_{\nu_d\in \mu_{n_d}}(R^1j_*\cM(i_0,\ldots,i_{d-1}))_{x_d=\nu_d}.
\]
Let $\cV_\HG(a_0,\ldots,a_{d-1})$ be the locally constant sheaf associated
to the monodromy representation \eqref{per-lem1-eq1} of
the hypergeometric function 
\[
{}_dF_{d-1}\left({a_0,\ldots,a_{d-1}\atop 1,\ldots,1};t\right).
\] 
Let $\rho_\alpha:S_\alpha\to \bA^1(t)\setminus\{0,1\}$ be the morphism such that
$\rho_\alpha^*(t)=\alpha/(1-x_d^{n_d})$.
Then it folllows from (iii)$_{d-1}$ and Corollary \ref{dim-cor} that
$\cM(i_0,\ldots,i_{d-1})$ is isomorphic to the dual of
$\rho^*\cV_\HG(a_0,\ldots,a_{d-1})$. 
Therefore
\[
(R^1j_*\cM(i_0,\ldots,i_{d-1}))_{x_d=\nu_d}=\Coker[T_\infty-\id:
V_\HG(a_0,\ldots,a_{d-1})_\alpha
\to V_\HG(a_0,\ldots,a_{d-1})_\alpha]
\]
where $T_p$ is the local monodromy at $t=p$.
As is well-known, the right hand side vanishes.
This completes the proof of (iii)$_d$.
\end{pf}
\begin{thm}\label{dim-thm3}
Let $n_k\geq 1$.
Let $Q$ be the set of $(d+1)$-tuple $(i_0,\ldots,i_d)$ of integers 
such that $0\leq i_k< n_k$ for all $k$ and $i_k=0$ for some $k$.
Then for all $0\leq j\leq d$,
the subspace
\begin{equation}\label{dim-thm3-eq1}
\bigoplus_{(i_0,\ldots,i_d)\in  Q}H^j_\dR(U_\alpha/\C)(i_0,\ldots,i_d)
\end{equation}
of the de Rham cohomology $H^j_\dR(U_\alpha/\C)$
is generated by exterior products of
\begin{equation}\label{dim-thm3-eq2}
\dlog(x_i-\nu_i)=\frac{dx_i}{x_i-\nu_i},\quad \nu_i\in \mu_{n_i}.
\end{equation}
More precisely, for $(i_0,\ldots,i_d)\in Q$, let us put 
$I=\{k\in \{0,1,\ldots,d\}\mid i_k>0\}$
and $I^c:=\{0,1,\ldots,d\}\setminus I$.
Put $s:=\sharp I$
and 
\[
\omega:=\bigwedge_{k\in I}\frac{x^{i_k-1}_kdx_k}{x_k^{n_k}-1}\in 
H^s_\dR(U_\alpha/\C)(i_0,\ldots,i_d)
\]
a $s$-form.
Let $p_l:U_\alpha\to T_l:=\bA^1(y_l)\setminus\{y_l=1\}$ be the morphism
given by $p_l^*(y_l)=x_l^{n_l}$. Fix $k_0\in I^c$.
Put
$T:=\prod_{l\in I^c\setminus\{k_0\}}T_l$, and
$p:=\prod p_l:U_\alpha\to T$.
Then the map
\begin{equation}\label{dim-thm3-eq3}
H^{j-s}_\dR(T/\C)\lra H^j_\dR(U_\alpha/\C)(i_0,\ldots,i_d),
\quad x\mapsto \omega\wedge p^*x
\end{equation}
is bijective. In particular
\[
\dim H^j_\dR(U_\alpha/\C)(i_0,\ldots,i_d)
=\binom{d-s}{j-s}=\binom{d-s}{d-j}.
\]
\end{thm}
\begin{pf}
Fix a primitive root $\nu_k\in \mu_{n_k}$ for each $k$. 
Let $\sigma_k:=\sigma_k(\nu_k)$ be the automorphism of $U_\alpha$ 
given by $(x_0,\ldots,x_d)\mapsto
(x_0,\ldots,\nu_kx_k,\ldots,x_d)$.
Let $U_\alpha/G$ denote the quotient scheme
by a finite group $G\subset \langle\sigma_0,\ldots,\sigma_d\rangle$.
Put $S_i:=\bA^1(x_i)\setminus\{x_i^{n_i}=1\}$.
Then we have
\begin{align*}
U_\alpha/\langle\sigma_k\rangle
&=\Spec\left[\C[x_0,\ldots,x_d]/
((1-x_0^{n_0})\cdots(1-x_d^{n_d})-\alpha)\right]^{\sigma=\id}\\
&\cong\Spec\C[x_0,\ldots,y_k,\ldots,x_d]/
((1-x_0^{n_0})\cdots(1-y_k)\cdots(1-x_d^{n_d})-\alpha)\\
&\cong\Spec\C[x_0,\ldots,\wh{x}_k,\ldots,x_d]
[(1-x_0^{n_0})^{-1}\cdots\wh{(1-x^{n_k}_k)}^{-1}\cdots(1-x_d^{n_d})^{-1}]\\
&\cong 
S_0\times\cdots\times\wh{S_k}\times\cdots \times S_d.
\end{align*}
This yields an isomorphism
\begin{equation}\label{dim-thm3-eq4}
H^j_\dR(U_\alpha/\C)^{\sigma_k=\id}\cong 
H^j_\dR(S_0\times\cdots\times\wh{S_k}\times\cdots \times S_d).
\end{equation}
The de Rham cohomology group
$H^1_\dR(S_i/\C)$ has a basis $\{dx_i/(x_i-\nu_i)\mid\nu_i\in \mu_{n_i}\}$.
Hence the subspace
\[
\bigoplus_{(i_0,\ldots,i_d)\in  Q}H^j_\dR(U_\alpha/\C)(i_0,\ldots,i_d)
=\sum_{k=0}^dH^j_\dR(U_\alpha/\C)^{\sigma_k=\id}
\]
is generated by exterior products of $dx_i/(x_i-\nu_i)$'s.
The eigen component $H^j_\dR(U_\alpha/\C)(i_0,\ldots,i_d)$ is the image of
a projector
\[
\prod_{0\leq k\leq d,\,k\ne k_0}\left(\frac{1}{n_k}\sum_{r=0}^{n_k-1}\nu^{-i_kr}_k\sigma_k^r\right)
\]
on $H^j_\dR(U_\alpha/\C)^{\sigma_{k_0}=\id}$.
Noticing the natural isomorphism \eqref{dim-thm3-eq4},
it is not hard to see that the image of the above projector is isomorphic to 
the left hand side of \eqref{dim-thm3-eq3}.
\end{pf}
\begin{thm}\label{dim-thm4}
Suppose $n_k>1$ and let $0<i_k<n_k$.
Let $F^\bullet$ be the Hodge filtration in the mixed Hodge structure
$H^d(U_\alpha,\Q)$.
Then $\dim\mathrm{Gr}_F^pW_dH^d(U_\alpha,\C)(i_0,\ldots,i_d)=1$ 
for each $p\in\{0,\ldots, d\}$.
\end{thm}
\begin{pf}
Let $f:U\to S$ be the fibration, and let $\cH:=W_dR^df_*\Q$ be the admissible variation of
Hodge structure of pure weight $d$.
We write $H_{B,\alpha}=\cH_\alpha=W_dH^d(U_\alpha,\Q)$.
Let $T$ be the local monodromy on $W_dH^d(U_\alpha,\Q)$ at $t=0$,
and put $N:=\log T$.
Let $W(N)$ be the monodromy weight filtration (cf. \cite[p.106--107]{morrison}).
Recall from Corollary \ref{dim-cor} that 
$H_{B,\alpha}(i_0,\ldots,i_d):=W_dH^d(U_\alpha,\C)(i_0,\ldots,i_d)$ is a $(d+1)$-dimensional representation of 
$\pi_1(S,\alpha)$ which is isomorphic to the dual of the
monodromy representation of the hypergeometric function
${}_{d+1}F_d\left({a_0,\ldots,a_d\atop 1,\ldots,1};t\right)$.
In particular, $T$ is unipotent 
and $N$ on $H_{B,\alpha}(i_0,\ldots,i_d)$ has rank $d$.
Therefore 
\begin{equation}\label{dim-thm4-eq1}
\dim\mathrm{Gr}^{W(N)}_{2p}H_{B,\alpha}(i_0,\ldots,i_d)=1
\end{equation}
for each $0\leq p\leq d$,
and the map
\begin{equation}\label{dim-thm4-eq2}
N:
\mathrm{Gr}^{W(N)}_{2p}H_{B,\alpha}(i_0,\ldots,i_d)\lra
\mathrm{Gr}^{W(N)}_{2p-2}H_{B,\alpha}(i_0,\ldots,i_d)
\end{equation}
is bijective for each $1\leq p\leq d$.
Let $\cH_{\lim}$ be the limiting mixed Hodge structure of $\cH$ at $t=0$
with the limiting Hodge filtration $F^\bullet_{\lim}$, whose
underlying $\Q$-module is $H_{B,\alpha}$ and the weight filtration is given by
$W(N)$ (\cite{Schmid}, \cite[\S 5]{morrison}).
Since \eqref{dim-thm4-eq2} is bijective, we have
\begin{equation}\label{dim-thm4-eq3}
\xymatrix{
\mathrm{Gr}^{W(N)}_{2d} H_{B,\alpha}\ar[r]^N_\sim
&\mathrm{Gr}^{W(N)}_{2d-2} H_{B,\alpha}
\ar[r]^{\hspace{0.7cm}N}_{\hspace{0.7cm}\sim}
&\cdots\ar[r]^{N\hspace{0.7cm}}_{\sim\hspace{0.7cm}}
&\mathrm{Gr}^{W(N)}_0 H_{B,\alpha}
}
\end{equation}
where all arrows are morphisms of mixed Hodge structures of type $(-1,-1)$.
Since the weight piece
$\mathrm{Gr}^{W(N)}_0 \cH_{\lim}$ is of Hodge type $(0,0)$, we have that 
$\mathrm{Gr}^{W(N)}_{2p} \cH_{\lim}$ is of Hodge type $(p,p)$ for each $0\leq p\leq d$.
Therefore, 
there is a natural bijection
\[
\mathrm{Gr}^p_{F_{\lim}}(H_{B,\alpha}\ot\C)
\os{\cong}{\lra} \mathrm{Gr}^{W(N)}_{2p}(H_{B,\alpha}\ot\C)
\]
and this induces
\begin{equation}\label{dim-thm4-eq4}
\mathrm{Gr}^p_{F_{\lim}}H_{B,\alpha}(i_0,\ldots,i_d)
\os{\cong}{\lra} \mathrm{Gr}^{W(N)}_{2p}H_{B,\alpha}(i_0,\ldots,i_d).
\end{equation}
We have
\[
\dim\mathrm{Gr}^p_{F}[W_dH^d(U_\alpha,\C)(i_0,\ldots,i_d)]=
\dim\mathrm{Gr}^p_{F_{\lim}}[W_dH^d(U_\alpha,\C)(i_0,\ldots,i_d)]=1
\]
by \eqref{dim-thm4-eq1} and \eqref{dim-thm4-eq4}.
\end{pf}
\begin{cor}
Put $D:=t\frac{d}{dt}$.
Let $\omega_{i_0\ldots i_d}\in \vg(X,\Omega^d_{X/A})$ be the regular $d$-form
\eqref{form-eq1}. 
Then for $0\leq p\leq d$, 
\[\mathrm{Gr}_F^pW_dH^d_\dR(U/A)(i_0,\ldots,i_d)\] is a free $A$-module
of rank one
with basis
$D^{d-p}\omega_{i_0\ldots i_d}$.
\end{cor}
\begin{pf}
By Corollary \ref{dim-cor2}, $W_dH^d_\dR(U/A)(i_0,\ldots,i_d)$ is a free $A$-module
with basis $\{D^k\omega_{i_0\ldots i_d}\mid k=0,\ldots, d\}$.
Therefore $\mathrm{Gr}_F^0W_dH^d_\dR(U/A)(i_0,\ldots,i_d)$ is generated by 
$D^d\omega_{i_0\ldots i_d}$.
Theorem \ref{dim-thm4} shows that $\mathrm{Gr}_F^0W_dH^d_\dR(U/A)(i_0,\ldots,i_d)$
is locally free of rank one, and hence this is free with basis $D^d\omega_{i_0\ldots i_d}$.
Then $F^1W_dH^d_\dR(U/A)(i_0,\ldots,i_d)$ is locally free of rank $d$
generated by $\{D^k\omega_{i_0\ldots i_d}\mid k=0,\ldots, d-1\}$, and
$\mathrm{Gr}_F^1W_dH^d_\dR(U/A)(i_0,\ldots,i_d)$ is generated by 
$D^{d-1}\omega_{i_0\ldots i_d}$.
Again by Theorem \ref{dim-thm4}, one concludes that this is free of rank one with basis
$D^{d-1}\omega_{i_0\ldots i_d}$. Keep continuing the same argument,
one has the desired assertion for each $p$.
\end{pf}
\begin{center}
{\bf\large Summary}
\end{center}
\begin{enumerate}
\renewcommand{\labelenumi}{\rm{(\arabic{enumi})}}
\item
Let $0\leq j<d$.
Then $H^j_\dR(U_\alpha/\C)$
is generated by exterior products of $dx_i/(x_i-\nu_i)$'s.
In particular $H^j_B(U_\alpha,\Q)$ carries a Tate-Hodge structure of type $(j,j)$.
\item
Let $P$ (resp. $Q$) be the set of $(d+1)$-tuple $(i_0,\ldots,i_d)$ of integers 
such that $0< i_k< n_k$ for all $k$ (resp. $0\leq i_k< n_k$ and $i_k=0$ for some $k$).
Then
\[
H^d_B(U_\alpha,\C)=
\overbrace{\bigoplus_{(i_0,\ldots,i_d)\in P}H^d_B(U_\alpha,\C)(i_0,\ldots,i_d)}^{\text{weight } d}
\op
\overbrace{\bigoplus_{(i_0,\ldots,i_d)\in Q}H^d_B(U_\alpha,\C)(i_0,\ldots,i_d).}^{\text{weight } 2d}
\]
The weight $2d$ piece
is generated by exterior products of $dx_i/(x_i-\nu_i)$'s.
For $(i_0,\ldots,i_d)\in P$, one has the Hodge decomposition
\[
H^d_B(U_\alpha,\C)(i_0,\ldots,i_d)=\bigoplus_{p=0}^dH^{p,d-p},\quad \dim H^{p,d-p}=1.
\]
\item
For $(i_0,\ldots,i_d)\in P$,
the $\pi_1(S,\alpha)$-representation $H^d_B(U_\alpha,\C)(i_0,\ldots,i_d)$ is isomorphic to
the dual of the monodromy representation \eqref{per-lem1-eq1} of
the hypergeometric function 
\[
{}_{d+1}F_{d}\left({a_0,\ldots,a_d\atop 1,\ldots,1};t\right),\quad a_k:=1-\frac{i_k}{n_k}.
\] 
\item
Let $\cD:=\C\langle t,(t-t^2)^{-1},\frac{d}{dt}\rangle$ be the Wyle algebra.
Put $D:=t\frac{d}{dt}$ and $P_\HG:=D^{d+1}-t(D+a_0)\cdots(D+a_d)$.
Let $\omega_{i_0\ldots i_d}\in \vg(X,\Omega^d_{X/A})$ be the regular $d$-form
\eqref{form-eq1}. 
Then
\[
\cD/\cD P_\HG\os{\cong}{\lra}
H^d_\dR(U/A)(i_0,\ldots,i_d),\quad P\longmapsto P(\omega_{i_0\ldots i_d}).
\]
\item
For $0\leq j<d$,
\[
\dim H^j_B(U_\alpha,\C)(i_0,\ldots,i_d)=\begin{cases}
0&(i_0,\ldots,i_d)\in P\\
\binom{d-s}{d-j}&(i_0,\ldots,i_d)\in Q
\end{cases}
\]
where $s=\sharp\{k\in\{0,1,\ldots,d\}\mid i_k>0\}$, and for $j=d$
\[
\dim H^d_B(U_\alpha,\C)(i_0,\ldots,i_d)=\begin{cases}
d+1&(i_0,\ldots,i_d)\in P\\
1&(i_0,\ldots,i_d)\in Q.
\end{cases}
\]
\end{enumerate}
\section{Higher Ross Symbols}\label{HRoss-sect}

\subsection{Definition of Higher Ross Symbols}\label{Ross-defn-sect}
Let $U$ be the hypergeometric scheme over $A$ as in \S \ref{setting-sect}. 
Let
$\xi_0:=\{1-x_0,1-x_1,\ldots,1-x_d\}\in K^M_{d+1}(\O(U))$
be a symbol in Milnor's $K$-group.
Let $\sigma_i(\nu_i)$ be the automorphisms in \eqref{setting-eq1}, which act 
on the Milnor $K$-group in the natural way.
For $n_i$-th roots $\nu_k\in \mu_{n_k}(A)$ of unity,
we define
\begin{align}\label{Ross-defn-eq2}
\Ross=\Ross(\nu_0,\ldots,\nu_d)&:=
(1-\sigma_0(\nu_0))\cdots(1-\sigma_d(\nu_d))(\xi_0)\\
&=\left\{\frac{1-x_0}{1-\nu_0x_0},\frac{1-x_1}{1-\nu_1x_1},\cdots,\frac{1-x_d}{1-\nu_dx_d}
\right\}\in K^M_{d+1}(\O(U))\label{Ross-defn-eq1}
\end{align}
and call it a {\it higher Ross symbol}.
There is the canonical map $\iota:K^M_{d+1}(\O(U))\to K_{d+1}(U)$ to Quillen's $K$-group.
We also denote the element $\iota(\Ross)$ by the same notation.

\begin{lem}\label{Ross-dlog}
Let $N=\mathrm{lcm}(n_0,\ldots,n_d)$ and $\zeta_N\in \ol\Q$ a $N$-th 
primitive root of unity. 
Let $A=\Z[1/N,\zeta_N][t,(t-t^2)^{-1}]$. Let
\[
\dlog:K_s^M(\O(U))\to \vg(U,\Omega^s_{U/\Z[1/N,\zeta_N]})
\]
be the dlog map.
Then
\[
\dlog(\Ross)=(-1)^d\sum_{i_0=1}^{n_0-1}\cdots
\sum_{i_d=1}^{n_d-1}(1-\nu_0^{i_0})\cdots(1-\nu_d^{i_d})
\omega_{i_0\ldots i_d}\frac{dt}{t}
\]
where $\omega_{i_0\ldots i_d}$ are the rational differential forms in \S \ref{form-sect}.
\end{lem}
\begin{pf}
Exercise (left to the reader).
\end{pf}
\begin{center}
{\bf Higher Ross vs. Ross.}
\end{center}
Let us see the relation with the original Ross symbols.
Let $F$ be the Fermat curve defined by an equation $z^n+w^m-1$.
Recall from \cite[Theorem 2]{ross2} the Ross symbol
\[
\{1-z,1-w\}
\]
in Milnor's $K_2^M$ of the function field of $F$.
It is not hard to see that this has no boundary, namely it lies in the image of
$K_2(F)\ot\Q$.
Let $U$ be the one-dimensional hypergeometric scheme with $t=1$ 
which is an affine scheme defined by an equation $(1-x_0^{n_0})(1-x_1^{n_1})=1$
$\Leftrightarrow$ $x_0^{-n_0}+x_1^{-n_1}=1$.
Put $z=x_0^{-1}$ and $w=x_1^{-1}$.
Then our higher Ross symbol is
\[
\Ross(\nu_0,\nu_1)
=\left\{\frac{1-x_0}{1-\nu_0x_0},\frac{1-x_1}{1-\nu_1x_1}\right\}
=\left\{\frac{z-1}{z-\nu_0},\frac{w-1}{w-\nu_1}\right\}.
\]
Take the summation over all $n_0$-th and $n_1$-th roots ($\nu_0,\nu_1)$
of unity. Then we have
\begin{align*}
\sum_{\nu_0,\nu_1}\Ross(\nu_0,\nu_1)
&=\left\{\frac{(z-1)^{n_0}}{z^{n_0}-1},\frac{(w-1)^{n_1}}{w^{n_1}-1}\right\}\\
&=\{(z-1)^{n_0},(w-1)^{n_1}\}-\{(z-1)^{n_0},w^{n_1}-1\}
-\left\{z^{n_0}-1,\frac{(w-1)^{n_1}}{w^{n_1}-1}\right\}\\
&=\{(z-1)^{n_0},(w-1)^{n_1}\}-\overbrace{\{(z-1)^{n_0},-z^{n_0}\}}^{\text{2-torsion}}
-\overbrace{\left\{-w^{n_1},\frac{(w-1)^{n_1}}{w^{n_1}-1}\right\}}^{\text{2-torsion}}\\
&=n_0n_1\{1-z,1-w\}
\end{align*}
the original Ross symbol, up to 2-torsion.
\begin{rem}[cf. {\cite[(4.18)]{otsubo-2}}]\label{Ross-rem}
In \cite{ross1}, Ross considered another symbol
\[
\{1-zw,z\}
\]
in $K_2$ of the Fermat curve $F$ defined by $z^n+w^n=1$. However this is essentially 
the above his symbols. More precisely, let $C$ be the smooth projective curve
defined by an equation $y^n=x(1-x)$.
Let $\rho:F\to C$ be the covering given by $\rho^*y=zw$ and $\rho^*x=z^n$.
Then using the algorithm in \cite{RT}, one can show 
\begin{equation}\label{Ross-rem-eq1}
\rho_*\{1-z,1-w\}=3\{1-y,x\}
\end{equation}
up to torsion, and hence
\begin{equation}\label{Ross-rem-eq2}
\{1-zw,z\}=\frac13\rho^*\rho_*\{1-z,1-w\}
=\frac13\sum_{\zeta\in\mu_n(\ol\Q)}\{1-\zeta z,1-\zeta^{-1} w\}.
\end{equation}
\end{rem}
\subsection{Boundary of higher Ross symbols}\label{boundary-sect}
Let $X\supset U$ be a smooth compactification (Proposition \ref{sc-thm}).
By a standard argument using de Jong's alteration, one has that
the image $K_s(X)\ot\Q\lra K_s(U)\ot\Q$
does not depend on the choice of $X$.
We say that $u\in K_s(U)\ot\Q$ has (non trivial) boundary if it does not lie in the image of
$K_s(X)\ot\Q$.

The symbol $\xi_0\in K_{d+1}(U)$ has boundary, while we expect that the higher Ross symbols $\Ross$ have no boundary.
Several results convince the author of the truth, though
he has not succeeded to prove it.
\begin{prop}\label{boundary-lem}
Let $U\subset X$ be the smooth compactification in Proposition \ref{sc-thm}
and $Z=X\setminus U$.
Let $D_k=p_k^{-1}(x_k=\infty)_\red$ be the reduced 
fiber at $x_k=\infty$ where $p_k:X\to \P^1(x_k)$ is the projection.
Let $D_k^c$ be the union of irreducible components of $Z$ which are not
contained in $D_k$. Put $D^\circ_k=D_k\setminus(D_k\cap D_k^c)$.
Let $\partial:K_{*+1}(U)\to K'_*(Z)$ denote the boundary map in $K'$-theory.
Then 
\begin{equation}\label{boundary-lem-eq1}
\partial(\Ross)|_{D^\circ_k}=0
\end{equation}
in $K'_d(D^\circ_k)$ for any $k$. In particular
\[
\partial(\Ross)|_{Z^{\mathrm{sm}}}=0
\]
in $K'_d(Z^{\mathrm{sm}})$ where
$Z^{\mathrm{sm}}\subset Z$ is the maximal locus
which is smooth over $A$, or equivalently $Z^{\mathrm{sm}}=Z\setminus\cup_{i<j}
Z_i\cap Z_j$.
\end{prop}
\begin{pf}
Let $V_k\supset D_k$ be a formal neighborhood. 
Let $i:D_k\to V_k$ be the closed immersion.
Then there is a commutative diagram
\[
\xymatrix{
K_{d+1}(U\cap V_k)\ar[r]^\partial& K'_d(D_k^\circ)\\
K_d(U)\ot \O(V_k)^\times\ar[r]_{\partial\ot i^*\quad}\ar[u]&
K'_{d-1}(D_k^\circ)\ot \O(D_k)^\times.\ar[u]
}
\]
Recall from the definition that $\Ross$ is the image of an element
\[
(-1)^{d-k}\left\{
\frac{1-x_0}{1-\nu_0x_0},\cdots,
\widehat{\frac{1-x_k}{1-\nu_kx_k}},\cdots,\frac{1-x_d}{1-\nu_dx_d}
\right\}\ot
\left\{\frac{1-x_k}{1-\nu_kx_k}\right\}\in 
K_{d+1}(U)\ot \O(V_k)^\times.
\]
Since
\[
i^*\left(\frac{1-x_k}{1-\nu_kx_k}\right)=
\left(\frac{1-x_k}{1-\nu_kx_k}\right)\bigg|_{x_k=\infty}=\nu_k^{-1}
\]
is a torsion, one has the vanishing \eqref{boundary-lem-eq1}.
\end{pf}
\begin{cor}\label{boundary-prop1}
When $d=1, 2$, the higher Ross symbol $\Ross$
has no boundary.
\end{cor}
\begin{pf}
We may assume $A=\Z[1/N,\zeta_N][t,(t-t^2)^{-1}]$
with $N:=\mathrm{lcm}(n_0,\ldots,n_d)$.
Write $\O:=\Z[1/N,\zeta_N]$.
It follows from the localization sequence
\[
\xymatrix{
0\ar[r]&K_q(\O[t])^{(j)}\ar[r]&K_q(A)^{(j)}\ar[r]&(K_{q-1}(\O)^{(j-1)})^{\op 2}
\ar[r]&0\\
&K_q(\O)^{(j)}\ar[u]^\cong
}
\]
and Borel's theorem $K_q(\O)^{(j)}=0$ for $2j-q\ne1$ that we have
\begin{equation}\label{boundary-eq2-1}
K_q(A)^{(j)}=\begin{cases}
K_q(\O)^{(j)}&2j-q=1\\
(K_{q-1}(\O)^{(j-1)})^{\op2}&2j-q=2\\
0
\end{cases}
\end{equation}
for $q\geq 1$.
We want to show the vanishing
\begin{equation}\label{boundary-eq1}
\partial(\Ross)=0\in K'_d(Z)^{(d)}.
\end{equation}
If $d=1$, then this is immediate from Proposition \ref{boundary-lem}.
When $d=2$, one has
\[
\partial(\Ross)\in \Image[K_2(T)^{(1)} \to K'_2(Z)^{(2)}], 
\]
where $T:=Z\setminus Z^{\mathrm{sm}}=\cup_{i<j}Z_i\cap Z_j$ is disjoint union of $\Spec A$.
Hence it vanishes by \eqref{boundary-eq2-1}.
\end{pf}

\begin{prop}\label{boundary-prop2}
Let $N$ be the l.c.m. of $n_0,\ldots,n_d$, and let $\zeta_N\in \ol\Q$ be a primitive
$N$-th root of unity. Let $A=\Z[1/N,\zeta_N][t,(t-t^2)^{-1}]$. 
Let $l$ be a prime which divides $N$.
Suppose that the higher Chern class map
\begin{equation}\label{boundary-prop2-eq1}
K_d^Z(X)\ot\Q\lra H^{d+2}_{\et,Z}(X,\Q_l(d+1))
\end{equation}
to the etale cohomology group is injective.
Then 
the higher Ross symbols $\Ross$ have no boundary.
\end{prop}
\begin{pf}
It is well-known that
the diagram
\[
\xymatrix{
K_{d+1}(U)\ar[r]^\partial\ar[d]&K_d^Z(X)\ar[d]\\
H^{d+1}_\et(U,\Q_l(d+1))\ar[r]&H^{d+2}_{\et,Z}(X,\Q_l(d+1))
}
\]
is commutative. 
Since we assume the injectivity of \eqref{boundary-prop2-eq1},
it is enough to show that $\Ross$ vanishies in $H^{d+2}_{\et,Z}(X,\Q_l(d+1))$.
Put $P:=(1-\sigma_0(\nu_0))\cdots(1-\sigma_d(\nu_d))$ as in \eqref{Ross-defn-eq2}.
Since $\Ross=P(\xi_0)$,
it is enough to show that the composition
\[
\xymatrix{
H^{d+1}_\et(U,\Q_l(d+1))\ar[r]^P&H^{d+1}_\et(U,\Q_l(d+1))\ar[r]&H^{d+2}_{\et,Z}(X,\Q_l(d+1))
}
\]
is zero.
Let $f:X\to \Spec A$, $j:U\to X$ and $i:Z\to X$.
It is enough to show that the composition
\[
R(fj)_*\Q_l\os{P}{\lra}R(fj)_*\Q_l\lra R(fi)_*i^!\Q_l[1]
\]
is zero in the derived category of complexes of $l$-adic sheaves on $\Spec A$,
which is equivalent to that
\[
R^k(fj)_*\Q_l\os{P}{\lra}R^k(fj)_*\Q_l\lra R^{k+1}(fi)_*i^!\Q_l
\]
is zero for all $k\geq 0$.
Since all terms are smooth sheaves, it is enough to see this at one geometric point.
Then, thanks to the comparison theorem with the Betti or de Rham cohomology,
the assertion is reduced to that the composition
\[
H^k_\dR(U_\alpha/\C)\os{P}{\lra}
H^k_\dR(U_\alpha/\C)\lra H^{k+1}_{\dR,Z}(X_\alpha/\C)
\]
is zero for the fiber $U_\alpha$ at a point $t=\alpha\in\C\setminus\{0,1\}$.
However, this is immediate from Theorem \ref{dim-thm2}.
Indeed it implies that $P=0$ if $k<d$ and
\[
\Image(P)=W_dH^d_\dR(U_\alpha/\C)=\Image[H^d_\dR(X_\alpha/\C)\to H^d_\dR(U_\alpha/\C)]
\]
if $k=d$, which agrees with the kernel of 
$H^d_\dR(U_\alpha/\C)\to H^{d+1}_{\dR,Z}(X_\alpha/\C)$.
\end{pf}

\section{Beilinson Regulators of Higher Ross Symbols}\label{Beilinson-sect}

\subsection{Complex analytic function $\cF_{a_1,\ldots,a_s}(t)$}\label{analytic-sect}
For an integer $s\geq 1$, and $a_1,\ldots,a_s\in\C\setminus\Z_{\leq0}$, we put
\[
F_{a_1,\ldots,a_s}(t):={}_{s+2}F_{s+1}\left({a_1+1,\ldots,a_s+1,1,1\atop 2,\ldots,2};t\right),
\]
\[
\cF_{a_1,\ldots,a_s}(t):=\sum_{k=1}^s(\psi(a_k)+\gamma)+\log(t)+a_1\cdots a_st 
F_{a_1,\ldots,a_s}(t)
\]
where $\gamma=-\Gamma'(1)$ is the Euler constant and
$\psi(t)=\Gamma'(t)/\Gamma(t)$ is the digamma function, cf. \cite[5.2.2]{NIST}.
By definition
\begin{equation}\label{conn-eq0}
t\frac{d}{dt}(\cF_{a_1,\ldots,a_s}(t))
={}_sF_{s-1}\left({a_1,\ldots,a_s\atop 1,\ldots,1};t\right).
\end{equation}
The following theorem is often refered to as the {\it connection formula}
which describes the behavior of $\cF_{a_1,\ldots,a_s}(t)$
around $t=\infty$.
This will play a key role in the proof of 
Thorem \ref{mainB-thm} (main theorem) below.

\begin{thm}\label{conn}
Suppose that $a_i\not\in \Z$ for any $i$ and $a_i-a_j\not\in\Z$ for any $i\ne j$
and $\mathrm{Re}(\sum_{i=1}^sa_i)<s$.
Put
\[
H_j(t):=a_j^{-1}(-t)^{a_j}{}_{s+1}F_{s}\left(
{a_j,\cdots,a_j\atop 1-a_1+a_j,\cdots,\widehat{1-a_j+a_j},\cdots,1-a_s+a_j,1+a_j};t\right),
\]
\[
C_j:=
(\Gamma(1-a_j))^{-s+1}
\prod_{1\leq k\leq s,\,k\ne j} \frac{\Gamma(a_k-a_j)}{\Gamma(a_k)}.
\]
Then we have
\begin{equation}\label{conn-eq00}
\cF_{a_1,\ldots,a_s}(t^{-1})=\pi i-\sum_{j=1}^sC_jH_j(t).
\end{equation}
\end{thm}
\begin{pf}
Let $z=t^{-1}$. Put
\[
G_j(z):=z\frac{d}{dz}H_j(z)
=(-z)^{a_j}{}_sF_{s-1}
\left({a_j,\cdots,a_j\atop 1-a_1+a_j,\cdots,\widehat{1-a_j+a_j},\cdots,1-a_s+a_j};z\right).
\]
Then 
\begin{equation}\label{conn-eq2}
{}_sF_{s-1}\left({a_1,\ldots,a_s\atop 1,\ldots,1};t\right)
=\sum_{j=1}^sC_jG_j(z)
\end{equation}
by \cite[16.8.8]{NIST}.
Therefore we have
\[
t\frac{d}{dt}(\cF_{a_1,\ldots,a_s}(t))
=z\frac{d}{dz}\left(\sum_{j=1}^sC_jH_j(z)\right)
\]
by \eqref{conn-eq0}, and this implies
\[
\log(t)+a_1\cdots a_stF_{a_1,\ldots,a_s}(t)+\sum_{j=1}^sC_jH_j(z)=\text{constant}.
\]
We compute the constant
\begin{equation}\label{conn-eq3}
C:=a_1\cdots a_sF_{a_1,\ldots,a_s}(1)+\sum_{j=1}^sC_jH_j(1)
\end{equation}
where the series $F_{a_1,\ldots,a_s}(t)$ and $H_j(t)$ are absolutely convergent
on $|t|=1$ by the 
assumption $\mathrm{Re}(\sum_{i=1}^sa_i)<s$.
Noticing $(c)_n=\Gamma(c+n)/\Gamma(c)$, one has
\begin{align}
\Gamma(c)\left[{}_{s+1}F_{s}\left({a_1,\cdots,a_s,c\atop 1,\cdots,1};1\right)-1\right]
&=\sum_{n=1}^\infty\frac{(a_1)_n\cdots(a_s)_n\Gamma(c+n)}{n!^{s+1}}\notag\\
&\to\sum_{n=1}^\infty\frac{(a_1)_n\cdots(a_s)_n\Gamma(n)}{n!^{s+1}}\quad (\text{as }c\to0)
\notag\\
&=\sum_{n=0}^\infty\frac{(a_1)_{n+1}\cdots(a_s)_{n+1}n!}{(n+1)!^{s+1}}
\notag\\
&=a_1\cdots a_s\cdot{}_{s+2}F_{s+1}\left({a_1+1,\cdots,a_s+1,1,1\atop 2,\cdots,2};1\right)
\notag\\
&=a_1\cdots a_sF_{a_1,\ldots,a_s}(1)\label{conn-eq1}.
\end{align}
Let
\[
G^*_j(z)
:=(-z)^{a_j}{}_{s+1}F_{s}
\left({a_j,\cdots,a_j\atop 1-a_1+a_j,\cdots,\widehat{1-a_j+a_j},\cdots,1-a_{s}+a_j,1-c+a_j};z\right),
\quad j\leq s\]
\[
G^*_{s+1}(z)
:=(-z)^{c}{}_{s+1}F_{s}
\left({c,\cdots,c\atop 1-a_1+c,\cdots,\cdots,1-a_{s}+c};z\right),
\]
\[
C^*_j:=
\prod_{k=1,k\ne j}^s \frac{\Gamma(a_k-a_j)}{\Gamma(a_k)}
\cdot(\Gamma(1-a_j))^{-s}\cdot \frac{\Gamma(c-a_j)}{\Gamma(c)}
,\quad j\leq s\]
\[
C^*_{s+1}:=
\prod_{k=1}^s \frac{\Gamma(a_k-c)}{\Gamma(a_k)}
\cdot(\Gamma(1-c))^{-s}.
\]
By \cite[16.8.8]{NIST}, one has
\begin{equation}\label{conn-eq4}
{}_{s+1}F_{s}\left({a_1,\ldots,a_s,c\atop 1,\ldots,1};t\right)
=\sum_{j=1}^{s+1}C^*_jG^*_j(z).
\end{equation}
By \eqref{conn-eq1} and 
\[
\lim_{c\to0}\Gamma(c)C^*_jG^*_j(1)=(-a_j)^{-1}
\prod_{k=1,k\ne j}^s \frac{\Gamma(a_k-a_j)}{\Gamma(a_k)}
\cdot(\Gamma(1-a_j))^{-s+1}a_jH_j(1)
=-C_jH_j(1)
\]
for $1\leq j\leq s$, one has
\[
\lim_{c\to0}\Gamma(c)\left[{}_{s+1}F_{s}\left({a_1,\ldots,a_s,c\atop 1,\ldots,1};1\right)-1
-\sum_{j=1}^{s}C^*_jG^*_j(1)\right]
=\overbrace{a_1\cdots a_sF_{a_1,\ldots,a_s}(1)+\sum_{j=1}^sC_jH_j(1)}^C
\]
where $C$ is the constant \eqref{conn-eq3}.
By \eqref{conn-eq4}, the left hand side is
\[
\lim_{c\to0}\Gamma(c)\left[C^*_{s+1}G^*_{s+1}(1)
-1\right]=
\lim_{c\to0}\Gamma(c)\left[e^{\pi ic}C^*_{s+1}
-1\right]
\]
where the equality follows from the fact that
$G^*_{s+1}(1)=e^{\pi i c}+O(c^{s+1})$.
Therefore, we have
\begin{align*}
C
&=\lim_{c\to0}\Gamma(c)\left[e^{\pi i c}\prod_{k=1}^s \frac{\Gamma(a_k-c)}{\Gamma(a_k)}
\cdot(\Gamma(1-c))^{-s}
-1\right]\\
&=\frac{d}{dc}\left[e^{\pi i c}\prod_{k=1}^s \frac{\Gamma(a_k-c)}{\Gamma(a_k)}
\cdot(\Gamma(1-c))^{-s}\right]\bigg|_{c=0}\\
&=\pi i-\sum_{k=1}^s\psi(a_k)+s\psi(1)
=\pi i-\sum_{k=1}^s(\psi(a_k)+\gamma)
\end{align*}
as desired.
\end{pf}
\begin{cor}\label{conn-cor}
Let $a_i\in\C\setminus\Z_{\leq0}$ satisfy
that $\mathrm{Re}(\sum_{i=1}^sa_i)<s$ (possibly $a_i-a_j\in\Z$).
Then around $t=0$, the analytic function
$\cF_{a_1,\ldots,a_s}(t^{-1})$ can be written as a $\C$-linear combination of
\[
(\log t)^it^{a_j}h_{ij}(t),\quad (i\geq0,\, j=1,2,\ldots,s)
\]
modulo $\Q(1)$ by some holomorphic functions $h_{ij}(t)$ at $t=0$.
Here ``$\log t$" does not appear when $a_i-a_j\not\in\Z$ for all $i\ne j$. 
\end{cor}
\begin{pf}
When $a_i-a_j\not\in\Z$ for all $i,j$, this is immediate from Theorem \ref{conn}.
In general, apply it for $a'_i:=a_i+\epsilon_i$ and take the limit $\epsilon_i\to0$.
Then ``$\log t$" appears as above (details are left to the reader).
\end{pf}
\begin{exmp}
We show
\[
\mathrm{Re}(\cF_{\frac12,\frac12}(t))=
\mathrm{Re}\left[
-\log 16+\log t+\frac{t}{4}{}_4F_3\left({\frac32,\frac32,1,1\atop 2,2,2};t\right)\right]
=\mathrm{Re}\left[
-2t^{-\frac12}{}_3F_2\left({\frac12,\frac12,\frac12\atop1,\frac32};t^{-1}\right)
\right]
\]
for $t\in\R_{>0}$.
Theorem \ref{conn} does not directly implies this as the assumption $a_i-a_j\not\in\Z$ 
is not satisfied.
Let $a_1=\frac12+\epsilon$ and $a_2=\frac12$ with $\epsilon>0$. Then 
we have from \eqref{conn-eq00} that
\[
\cF_{\frac12+\epsilon,\frac12}(t^{-1})=\pi i-(C_1H_1(t)+C_2H_2(t))
\]
where 
\[
H_1(t):=\left(\frac12+\epsilon\right)^{-1}(-t)^{\frac12+\epsilon}{}_3F_2\left(
{\frac12+\epsilon,\frac12+\epsilon,\frac12+\epsilon\atop 1+\epsilon,\frac32+\epsilon};t\right),\quad
H_2(t):=2(-t)^{\frac12}{}_3F_2\left(
{\frac12,\frac12,\frac12\atop 1-\epsilon,\frac32};t\right),
\]
\[
C_1:=
\frac{\Gamma(-\epsilon)}{\Gamma(\frac12)\Gamma(\frac12-\epsilon)},\quad
C_2:=
\frac{\Gamma(\epsilon)}{\Gamma(\frac12)\Gamma(\frac12+\epsilon)}.
\]
Taking the limit $\epsilon\to0$,
we have
\[
\cF_{\frac12,\frac12}(t^{-1})=\lim_{\epsilon\to0}
\cF_{\frac12+\epsilon,\frac12}(t^{-1})=\pi i
+2\pi^{-1}\log(-t)(-t)^{\frac12}{}_3F_2\left(
{\frac12,\frac12,\frac12\atop 1-\epsilon,\frac32};t\right)+(-t)^{\frac12}h(t)
\]
with $h(t)$ a holomorphic function at $t=0$ whose Maclaurin expansion has coefficients in $\R$.
Let $t>0$ and take the real part. Then this turns out
\[
\mathrm{Re}(\cF_{\frac12,\frac12}(t^{-1}))=
\mathrm{Re}\left[-2t^{\frac12}{}_3F_2\left(
{\frac12,\frac12,\frac12\atop 1-\epsilon,\frac32};t\right)\right].
\]
\end{exmp}
\subsection{Beilinson regulators of Higher Ross symbols}\label{mainB-sect}
Let $A=\C[t,(t-t^2)^{-1}]$ and put $S=\Spec A=\P^1(t)\setminus\{0,1,\infty\}$.
Let $f:U\to S$ be the hypergeometric scheme in \S \ref{setting-sect}.
Let
\begin{equation}\label{mainB-hc}
c_{i,j}:K_i(U)^{(j)}\lra \Ext^{2j-i}_{\MHM(U)} (\Q,\Q(j)),\quad i,j\geq0
\end{equation}
be the {\it Beilinson regulator maps (higher Chern class maps)}
to the extension groups of mixed Hodge modules by Saito \cite{msaito1}.
We only discuss the case $(i,j)=(d+1,d+1)$.
There is the natural map
\[
\Ext^{d+1}_{\MHM(U)} (\Q,\Q(d+1))\lra \Ext^1_{\VMHS(S)} (\Q,R^df_*\Q(d+1))
\]
as $R^jf_*\Q(d+1)=0$ for $j>d$ where
$\VMHS(S)$ denotes the category 
of admissible variations of mixed Hodge structures on $S$.
Let
\begin{equation}
\rho_B:K_{d+1}(U)\lra \Ext^1_{\VMHS(S)} (\Q,R^df_*\Q(d+1))
\end{equation}
be the composition with $c_{d+1,d+1}$, which is compatible
under the action of $G=\mu_{n_0}\times\cdots\mu_{n_d}$.
Fix a primitive $n_i$-th root of unity $\ve_i$ for each $i$.
Let 
\[
P:=\prod_{i=0}^d(1-\sigma_i(\ve_i))\in \Q[G].
\]
where $\sigma_i(\ve_i)=(1,\ldots,\ve_i,\ldots,1)\in G$ 
(cf. the notation in \S \ref{setting-sect}).
It follows from Theorem \ref{dim-thm2} (iii)
that $P$ gives a splitting of the inclusion $W_dR^df_*\Q\to R^df_*\Q$ and hence
one has 
a decomposition
\[
R^df_*\Q
=W_dR^df_*\Q\op \mathrm{Gr}^W_{2d}R^df_*\Q
\cong W_dR^df_*\Q\op \Q(-d)^\op
\]
by Theorem \ref{dim-thm3} (see also 
\S \ref{dim-thm2} Summary (2)).
Let 
\begin{equation}\label{mainB-eq1}
\Ross=\left\{\frac{1-x_0}{1-\nu_0x_0},\cdots,\frac{1-x_d}{1-\nu_dx_d}
\right\}\in K^M_{d+1}(\O(U))
\end{equation}
 be a higher Ross symbol introduced in \S \ref{Ross-defn-sect}.
We claim
\begin{equation}\label{mainB-lem0}
\rho_B(\Ross)\in \Ext^1_{\VMHS(S)}(\Q,W_dR^df_*\Q(d+1)).
\end{equation}
Indeed, 
let $\xi_0:=\{1-x_0,\ldots,1-x_d\}\in K^M_{d+1}(\O(U))$
and $P_{\ul \nu}:=(1-\sigma_0(\nu_0))\cdots(1-\sigma_d(\nu_d))$.
Then \eqref{mainB-lem0} follows from the fact that $\Ross=P_{\ul \nu}\xi_0$
and the compatibility of $\rho_B$ under the action of $G$ on noticing 
$\Image P_{\ul\nu}\subset\Image P$.

The element $\rho_B(\Ross)$ defines a 1-extension
\begin{equation}\label{mainB-ext}
\xymatrix{
0\ar[r]&W_dR^df_*\Q(d+1)\ar[r]&M_\Ross(U/S)\ar[r]&\Q_S\ar[r]&0.
}
\end{equation}
Let $U_\alpha=f^{-1}(\alpha)$ be the fiber at $t=\alpha$.
Restricting the 1-extension \eqref{mainB-ext} to $t=\alpha$, we have a 1-extension
\begin{equation}\label{mainB-ext-a}
\xymatrix{
0\ar[r]&W_dH^d(U_\alpha,\Q(d+1))\ar[r]&M_\Ross(U_\alpha)\ar[r]&\Q\ar[r]&0
}
\end{equation}
of mixed Hodge structures.
There are the natural isomorphism
\begin{align*}
\Ext^1_\MHS(\Q,W_dH^d(U_\alpha,\Q(d+1)))
&\cong W_dH^d(U_\alpha,\C/\Q(d+1))\\
&\cong \Hom(H_d(U_\alpha,\Q)/W_{-d-1},\C/\Q(d+1)).
\end{align*}
The extension class $\rho_B(\Ross)|_{U_\alpha}=\rho_B(\Ross|_{U_\alpha})$ defines an element
\begin{equation}\label{mainB-ext-d}
\reg_B(\Ross|_{U_\alpha})\in W_dH^d(U_\alpha,\C/\Q(d+1))
\end{equation}
in the Betti cohomology group.
In more explicit manner, this is defined in the following way.
Let $M_\Ross(U_\alpha)_B$ be the underlying $\Q$-module of 
$M_\Ross(U_\alpha)$, and $(M_\Ross(U_\alpha)_\dR,F^\bullet)$ the underlying
$\C$-module with the Hodge filtration.
There are two liftings of $1\in\Q$ in the exact sequence \eqref{mainB-ext-a}, say
$e_{B,\alpha}\in M_\Ross(U_\alpha)_B$ and $e_{\dR,\alpha}\in F^0M_\Ross(U_\alpha)_\dR$.
Then 
\[
\reg_B(\Ross|_{U_\alpha})=e_{\dR,\alpha}-e_{B,\alpha}\mod W_dH^d(U_\alpha,\Q(d+1)).
\]
\begin{thm}\label{mainB-thm}
For $0<i_k<n_k$ we write $a_k=1-i_k/n_k$.
Let $\alpha\in \C\setminus\{0,1\}$ and let $\Delta_\alpha\in H_d(U_\alpha,\Z)$ 
be the homology cycle in Theorem \ref{per-thm1}.
Then
\[
\langle\reg_B(\Ross|_{U_\alpha})\mid\Delta_\alpha\rangle=
\pm\sum_{0< i_k<n_k}(1-\nu_0^{i_0})\cdots(1-\nu_d^{i_d})
\frac{(2\pi i)^d}{n_0\cdots n_d}\cF_{a_0\ldots a_d}(\alpha)
\]
modulo $\Q(d+1)$, where $\langle-\mid-\rangle$ is the natural pairing
\[
H^d(U_\alpha,\C/\Q(d+1))\ot H_d(U_\alpha,\Q)\lra \C/\Q(d+1).
\]
\end{thm}
\begin{pf}
We think of $\langle\reg_B(\Ross|_{U_\alpha})\mid\Delta_\alpha\rangle$ to be a function on 
$\alpha$, which 
we write by $\langle\reg_B(\Ross)\mid\Delta_t\rangle$ with variable $t$.
This is a multi-valued function which
is locally holomorphic on $\C\setminus\{0,1\}$.
Recall from Lemma \ref{Ross-dlog} that
\[
\dlog(\Ross)=(-1)^d\sum_{0< i_k<n_k}(1-\nu_0^{i_0})\cdots(1-\nu_d^{i_d})
\omega_{i_0\ldots i_d}\frac{dt}{t}.
\]
It follows from Theorem \ref{per-thm1} together with \cite[Prop.3.1]{A}
that one has
\begin{align*}
&t\frac{d}{dt}\langle\reg_B(\Ross)\mid\Delta_t\rangle=
\sum_{0< i_k<n_k}(1-\nu_0^{i_0})\cdots(1-\nu_d^{i_d})
\int_{\Delta_t}\omega_{i_0\ldots i_d}\\
=&\pm\sum_{0< i_k<n_k}(1-\nu_0^{i_0})\cdots(1-\nu_d^{i_d})
\frac{(2\pi i)^d}{n_0\cdots n_d}\cdot
{}_{d+1}F_d\left({a_0,\ldots,a_d\atop1,\ldots,1};t\right).
\end{align*}
This implies
\begin{equation}\label{mainB-thm-eq1}
\langle\reg_B(\Ross)\mid\Delta_t\rangle=C+
\sum_{0< i_k<n_k}(1-\nu_0^{i_0})\cdots(1-\nu_d^{i_d})
\frac{(2\pi i)^d}{n_0\cdots n_d}\cF_{a_0\ldots a_d}(t)
\end{equation}
with a constant $C$.
There remains to show $C\in \Q(d+1)$. 
Let $T_\infty$ be the local momodromy $T_\infty$
at $t=\infty$.
We look at the action of
\[
Q:=\prod_{k=0}^d\prod_{i_k=1}^{n_k-1}(T_\infty-e^{2\pi i(1-i_k/n_k)})\in\Z[\pi_1(S,\alpha)].
\]
It follows from Corollary \ref{conn-cor} that $Q^m$ annihilates
all $\cF_{a_0\ldots a_d}(t)$ for sufficiently large $m>0$.
Since $T_\infty$ acts on the constant terms as identity, we have
\begin{equation}\label{mainB-thm-eq2}
Q^m\langle\reg_B(\Ross)\mid\Delta_t\rangle=Q^mC
=\left(\prod_{k=0}^d\prod_{i_k=1}^{n_k-1}(1-e^{2\pi i(1-i_k/n_k)})\right)^mC=(n_0\cdots n_d)^m C.
\end{equation}
On the other hand, since 
$\langle\reg_B(\Ross|_{U_\alpha})\mid\Delta_\alpha\rangle$ is defined for
$\Delta_\alpha$ modulo the kernel of $H_d(U_\alpha,\Q)\to H_d(X_\alpha,\Q)$,
it follows from Theorem \ref{dim-thm2} that one can replace $\Delta_\alpha$
with $\Delta'_\alpha$ where 
$\Delta_\alpha=\Delta'_\alpha+\Delta^{\prime\prime}_\alpha$
is the decomposition such that $\Delta^{\prime}_\alpha\in
\oplus_{0<i_k<n_k}H_d(U_\alpha,\C)(i_0,\ldots,i_k)$ and
$\Delta^{\prime\prime}_\alpha\in
\oplus_{\exists\, i_k=0}H_d(U_\alpha,\C)(i_0,\ldots,i_k)$.
We have
\[
Q^m\langle\reg_B(\Ross|_{U_\alpha})\mid\Delta_\alpha\rangle
=Q^m\langle\reg_B(\Ross|_{U_\alpha})\mid\Delta'_\alpha\rangle
=Q^m\cdot\overbrace{\langle e_{B,\alpha}\mid \Delta'_\alpha\rangle}^{(2\pi i)^{d+1}\Q}
-\langle e_{\dR,\alpha}\mid Q^m\Delta'_\alpha\rangle.
\]
For $0<i_k<n_k$, the component $H_d(U_\alpha,\C)(i_0,\ldots,i_k)$ is isomorphic
to the monodromy representation of the hypergeometric function by 
Corollary \ref{dim-cor}.
In particular it is annihilated by $Q^m$, and hence
we have $\langle e_{\dR,\alpha}\mid Q^m\Delta'_\alpha\rangle=0$.
Therefore the most left term in \eqref{mainB-thm-eq2} is contained in $\Q(d+1)$.
This completes the proof.
\end{pf}
\begin{cor}
Let $\ul\ve=(\ve_0,\ldots,\ve_d)\in G$, and put
$\Delta_\alpha(\ul\ve):=\sigma_0(\ve_0)\cdots\sigma_d(\ve_d)(\Delta_\alpha)$.
Then
\[
\langle\reg_B(\Ross|_{U_\alpha})\mid\Delta_\alpha(\ul\ve)\rangle=
\pm\sum_{0< i_k<n_k}\ve_0^{i_0}\cdots\ve_d^{i_d}(1-\nu_0^{i_0})\cdots(1-\nu_d^{i_d})
\frac{(2\pi i)^d}{n_0\cdots n_d}\cF_{a_0\ldots a_d}(\alpha)
\]
modulo $\Q(d+1)$.
\end{cor}
\begin{pf}
By the naturalness of the regulator map, we have
\[
\langle\reg_B(\Ross|_{U_\alpha})\mid\Delta_\alpha(\ul\ve)\rangle=
\langle\reg_B(\sigma_0(\ve_0)\cdots\sigma_d(\ve_d)
\Ross|_{U_\alpha})\mid\Delta_\alpha\rangle.
\]
We have
\begin{align*}
\sigma_0(\ve_0)\cdots\sigma_d(\ve_d)\Ross
&=
\left(\prod_{k=0}^d\sigma_k(\ve_k)-\sigma_k(\ve_k\nu_k)\right)
\xi_0\quad (\text{see }\eqref{Ross-defn-eq2})
\\
&=
\left(\prod_{k=0}^d1-\sigma_k(\ve_k)-(1-\sigma_k(\ve_k\nu_k))\right)\xi_0\\
&=
\sum(-1)^k \Ross(\ve_0\nu_0,\ldots,\ve_{i_1},\ldots,\ve_{i_k},\ldots,\ve_d\nu_d).
\end{align*}
The corollary follows by applying Theorem \ref{mainB-thm} to the Ross symbols 
in the last term.
\end{pf}
\section{Application to the Beilinson conjecture}\label{application-sect}
\subsection{Beilinson conjecture}
For a separated scheme $V$ of finite type over $\R$, 
let $H_B^\bullet(V,\Q)=H_B^\bullet(V(\C),\Q)$ denote 
the Betti cohomology groups. 
We denote by $F_\infty$ the infinite Frobenius on $V(\C)$ which is induced by
the complex conjugation.
It acts on
the Betti cohomology $H^\bullet_B(V,\Z)=H_B^\bullet(V(\C),\Z)$, and it naturally extends
on $H^\bullet_B(V,\C)$ as an anti-linear involution.
Moreover $F_\infty$ acts on the de Rham cohomology groups $H^\bullet_\dR(V/\R)\ot\C$
as an anti-linear involution such that it is identity on $H^\bullet_\dR(V/\R)$.
The actions of $F_\infty$ are compatible under the comparison 
$H^\bullet_\dR(V/\R)\ot\C\cong H^\bullet_B(V,\C)$.

Let $X$ be a smooth variety over $\R$. Let $A=\Z,\Q$ or $\R$.
Let $A(r)_\cD$ be the Deligne-Beilinson complex
\[
(2\pi\sqrt{-1})^rA_X\lra\O_X\lra\cdots\lra\Omega^{r-1}_X
\]
of sheaves on the analytic site $X^{an}$. 
The Deligne-Beilinson cohomology group is defined to be 
the hypercohomology group\[
H^\bullet_\cD(X,A(r)):={\mathbb H}^\bullet(X^{an},A(r)_\cD).
\]
We refer \cite{ev} or \cite[\S 2]{schneider}
for a general theory on the Deligne-Beilinson cohomology.
There is the natural map on $F_\infty^{-1}A(r)_\cD\to A(r)_\cD$
given by $\omega\mapsto\ol{F^*_\infty\omega}$, which we also denote by $F_\infty$
as long as there is no fear of confusion.
Let
\[
c^\cD_{i,j}:K_i(X)^{(j)}\lra H^{2j-i}_\cD(X,\R(j))^{F_\infty=1}
\]
be the Beilinson regulator map (or higher Chern class map)
to the fixed part of Deligne-Beilinson cohomology 
by $F_\infty$. This is compatible with the 
regulator map $c_{i,j}$ in \eqref{mainB-hc}
under the comparison
\[
\Ext^p_{\MHM(X)}(\Q,\Q(j))\cong H^p_\cD(X,\Q(j)).
\]
We now focus on the case $i=j>0$.
Then there is the natural isomorphism (e.g. \cite[p.9]{schneider})
\[
H^{j}_\cD(X,\R(j))^{F_\infty=1}\cong H_B^{j-1}(X,\R(j-1))^{F_\infty=1}
(=H_B^{j-1}(X,\R)^{F_\infty=(-1)^{j-1}}\ot\R(j-1)).
\]
The right hand side is endowed with the canonical $\Q$-structure
$H_B^{j-1}(X,\Q(j-1))^{F_\infty=1}$.
Composing the above with $c^\cD_{j,j}$, we have
\[
\reg_\R:K_j(X)^{(j)}\lra H^{j-1}_B(X,\R(j-1))^{F_\infty=1}, \quad j>0.
\]
\begin{conj}[Beilinson conjecture, \cite{beilinson}, \cite{schneider}]
\label{d=2-Bconj}
Let $X_\Q$ be a smooth projective variety over $\Q$, and $X_\R:=X_\Q\times_\Q\R$.
Let $K_i(X_\Q)^{(j)}_\Z\subset K_i(X_\Q)^{(j)}$ be the integral part
(\cite{Scholl}). Let $j>0$ be an integer.
Then 
\[
\reg_\R:K_j(X_\Q)^{(j)}_\Z\ot\R\lra H^{j-1}_B(X_\R,\R(j-1))^{F_\infty=1}
\]
is bijective, and 
\begin{equation}\label{d=2-Bconj-eq1}
\det[K_j(X_\Q)^{(j)}_\Z]\ot\det
[H^{j-1}_B(X_\R,\Q(j-1))^{F_\infty=1}]^{-1}\sim_{\Q^\times}L^*(h^{j-1}(X_\Q),0)
\end{equation}
where $L(M,0)$ denotes the $L$-function of a motive $M$, and 
where $L^*(M,m)$ is the leading coefficient in a Taylor series expansion at $s=m$.
Note that $\ord_{s=0} L(h^{j-1}(X_\Q),s)=\dim H^{j-1}_B(X_\R,\Q(j-1))^{F_\infty=1}$
under the hypothesis of the functional equation (\cite[p.5, Corollary]{schneider}).
\end{conj}
The left hand side of \eqref{d=2-Bconj-eq1} is called the {\it Beilinson regulator}, which
is a higher dimensional generalization of the classical Dirichlet regulators.
In more down-to-earth manner, it is described in the following way.
Let $\{\xi_p\}_p$ be a $\Q$-basis of $K_j(X_\Q)_\Z^{(j)}$ and
$\{\gamma_q\}_q$ a $\Q$-basis of $H_2^B(X,\Q)^{F_\infty=(-1)^{j-1}}$.
Let $\langle-\mid-\rangle$ denote the natural pairing on $H^\bullet_B(X)\ot H^B_\bullet(X)$.
Then
\begin{equation}\label{d=2-Bconj-eq2}
\text{LHS of }\eqref{d=2-Bconj-eq1}=
\det\left[\frac{1}{(2\pi \sqrt{-1})^{j-1}}
\langle\reg_\R(\xi_p)\mid\gamma_q\rangle\right]_{p,q}
\end{equation}
\subsection{$K_2$ of Elliptic curves}\label{EC-sect}
Let 
\[
U_\alpha=\Spec \Q[x_0,x_1]/((1-x_0^2)(1-x_1^2)-\alpha),\quad \alpha\in\Q\setminus\{0,1\}
\]
be a hypergeometric scheme in \S \ref{setting-sect}.
Let $X_\alpha$ be the smooth compactification of $U_\alpha$, which 
is an elliptic curve over $\Q$.
Then $\dim H^1_B(X_\alpha,\Q(1))^{F_\infty=1}=1$, and the Beilinson conjecture 
predicts that
there is an integral element $\xi\in K_2(X_\alpha)^{(2)}_\Z$ satisfying
\begin{equation}\label{EC-eq1}
\frac{1}{2\pi\sqrt{-1}}\langle\reg_\R(\xi)\mid \gamma\rangle\sim_{\Q^\times}L'(X_\alpha,0)
,\quad L(X_\alpha,s):=L(h^1(X_\alpha),s)
\end{equation}
where $\gamma\in H_1^B(X_\alpha,\Q)^{F_\infty=-1}$ is a generator 
(this statement is often referred to as the weak Beilinson conjecture).
Let $\Delta_\alpha$ be the
homology cycle in Theorem \ref{per-thm1}.
We take $\gamma:=\frac12(F_\infty+1)\Delta_\alpha$.
Let
\[
\Ross=\left\{\frac{1-x_0}{1+x_0},\frac{1-x_1}{1+x_1}\right\}
\]
be the higher Ross symbol.
Then it follows from Theorem \ref{mainB-thm} that one has
\[
\frac{1}{2\pi\sqrt{-1}}\langle\reg_\R(\Ross)\mid \gamma\rangle
=\mathrm{Re}[\cF_{\frac12,\frac12}(\alpha)].
\]
Thus we arrive at the following statement,
\begin{conj}
Suppose that $\Ross\in K_2(X_\alpha)^{(2)}$ is integral. Then
\[
\mathrm{Re}[\cF_{\frac12,\frac12}(\alpha)]\sim_{\Q^\times}L'(X_\alpha,0).
\]
\end{conj} 
If the denominator of $j(X_\alpha)=16(\alpha^2-16\alpha+16)^3/((1-\alpha)\alpha^4)$
is prime to $\alpha$ (e.g. $\alpha=\pm 2^n$, $n\in \{\pm 1,\pm 2,\pm 3\}$),
then $\Ross$ is integral.

\bigskip
\begin{center}
{\bf Numerical verifications}
\end{center}
 \begin{tabular}{c|c|c|c}
 $\alpha$&$\mathrm{Re}[\cF_{\frac12,\frac12}(\alpha)]$&$L'(X_\alpha,0)$&
 $\mathrm{Re}[\cF_{\frac12,\frac12}(\alpha)]/L'(X_\alpha,0)$\\
 \hline
$2$&$  -1.4866664931$&$0.74333324664$&$-2$\\
$-2$&$ -2.42449751304$ &$ 2.42449751304$&$-1$\\
$1/2$&$-3.3173289967$&
$  1.6586644983$&$-2$\\
$-1/2$&$  -3.5763399863$&$ -3.5763399863$&$1$\\
$4$&$  -1.0228481341$&$ 0.51142406705$&$-2$\\
$-4$&$  -1.942820350 $ &$0.971410175 $&$-2$\\
$1/4$&$ -4.091392536 $&$0.51142406705  $&$-8$\\
$-1/4$&$ -4.21743424174$&$2.10871712$&$-2$\\
$8$&$  -.71480404895$&$  1.429608097$&$-1/2$\\
$-8$&$  -1.5342722011 $ &$0.511424067 $&$-3$\\
$1/8$&$-4.819613084$&$ -9.639226168  $&$1/2$\\
$-1/8$&$ -4.8822409859$&$4.8822409859$&$-1$\\
\end{tabular}


\subsection{$K_3$ of K3 surfaces}\label{K3-sect}
In \cite{ono}, Ahlgren, Ono and Penniston study the $L$-function of 
a K3 surface $Y_\alpha$ over $\Q$ defined by an affine equation
\begin{equation}\label{ono-eq}
w^2=u_1u_2(1+u_1)(1+u_2)(u_1-\alpha u_2),\quad \alpha\in \Q\setminus\{0,1\}.
\end{equation}
The K3 surface $Y_\alpha$ is related to our hypergeometric scheme
\begin{equation}\label{K3-eq1}
U_\alpha:=\Spec \Q[x_0,x_1,x_2]/((1-x_0^2)(1-x_1^2)(1-x_2^2)-\alpha).
\end{equation}
To see this, we begin with a lemma in a general situation. 

\begin{lem}\label{K3-lem}
Let $A$ be a commutative ring.
Let $N\geq2$ be an integer, and let
\[
U_t:=\Spec A[x_0,\ldots,x_d]/((1-x^N_0)\cdots (1-x^N_d)-t)
\]
be a hypergeometric scheme with $t(1-t)\in A^\times$.
Let $n,n_1,\ldots,n_d$ be integers such that $0<n,n_i<N$ and
$\gcd(N,n,n_1,\ldots,n_d)=1$ for all $i$.
Let
\[
V_t:=\Spec A[y,z_1,\ldots,z_d]/(y^N-z_1^n\cdots z_d^n(1-z_1)^{n_1}\cdots 
(1-z_d)^{n_d}(-t+z_1\cdots z_d)^{N-n})
\]
be an affine scheme.
Then there is a covering morphism
\begin{equation}\label{covering}
\rho:U_t\lra V_t
\end{equation}
given by
\[
\begin{cases}
\rho^*y=x_0^{N-n}x_1^{n_1}\cdots x^{n_d}_d(1-x^N_1)\cdots (1-x^N_d)\\
\rho^*z_i=1-x_i^N\,(1\leq i\leq d).
\end{cases}
\]
\end{lem}
\begin{pf}
Straightforward.
\end{pf}

Let us take $A=\Q$, $t=\alpha$ and
$d=N=2$ and $n=n_1=n_2=1$ in Lemma \ref{K3-lem}.
Changing the variables
$u_1=-z_1$, $u_2=-1/z_2$ and $w=y/z_2^2$, we see that the equation of
$V_\alpha$ turns out to be
the equation \eqref{ono-eq}.
Therefore the smooth compactification of $V_\alpha$ is the K3 surface $Y_\alpha$
by Ahlgren, Ono and Penniston.
Let $X_\alpha$ be the smooth compactification of $U_\alpha$, which is a K3 surface
over $\Q$. 
The covering map \eqref{covering} induces a dominant rational morphism
\begin{equation}\label{covering2}
X_\alpha\lra Y_\alpha
\end{equation}
of K3 surfaces over $\Q$.

\medskip

For a smooth projective variety $S$ over $\Q$, we denote by $\NS(S)$ the
Neron-Severi group of $S\times_\Q\ol\Q$.
Let $h^2_\tr(S,\Q(m))=h^2(S,\Q(m))/\NS(S)\ot\Q(m-1)$ be the transcendental part of
the motive $h^2(S,\Q(m))$ (cf. \cite[7.2.2]{BMP}).
We denote by
$L(h^2_\tr(S),s)$ the $L$-function of $h^2_\tr(S,\Q)$.

\medskip

The rational morphism \eqref{covering2} induces an isomorphism
\begin{equation}\label{XY}
h^2_\tr(Y_\alpha,\Q)\cong
h^2_\tr(X_\alpha,\Q)
\end{equation}
of motives over $\Q$.
\begin{thm}[{\cite[Theorem 1.1]{ono}}]\label{d=2-thm1}
Let $\alpha\in\Q\setminus\{0,1\}$.
Let
\[
E_\alpha:y^2=x\left(x^2+2x-\frac{\alpha}{1-\alpha}\right)
\]
\[
E'_\alpha:(1-\alpha)y^2
=x\left(x^2+2x-\frac{\alpha}{1-\alpha}\right)
\]
be elliptic curves over $\Q$.
Then $L(h^2_\tr(Y_\alpha),s)=L(h^2_\tr(E_\alpha\times E'_\alpha),s)$.
\end{thm}
Van Geemen and Top construct an explicit correspondence between the $\Q$-motives
$h^2_\tr(Y_\alpha)$ and $h^2_\tr(E_\alpha\times E'_\alpha)$.
Let $Y_\alpha\to \P^1(u_2)$ be the morphism given by the projection 
$(w,u_1,u_2)\mapsto u_2$, which gives an elliptic fibration.
Let $N(Y_\alpha)$ be the $\Q$-linear subspace of $\NS(Y_\alpha)\ot\Q$ generated by components of
singular fibers and a section. Then the rank of $N(Y_\alpha)$ is $19$.
Let $\Delta'\subset E_\alpha\times E'_\alpha$ be the graph of
a morphism $(x,y)\mapsto (x,y/\sqrt{1-\alpha})$, and $\Delta'(1,1)\in \NS(
E_\alpha\times E'_\alpha)$
the K\"unneth $(1,1)$-component.
\begin{thm}[{\cite[Theorem 1.2]{GT}}]\label{GT}
There is an explicit correspondence between
$Y_\alpha$
and $E_\alpha\times E'_\alpha$ defined over $\Q$, which gives 
an isomorphism 
\begin{equation}
h^2(Y_\alpha)/N(Y_\alpha)\cong
h^1(E_\alpha)\ot h^1(E'_\alpha)/\Q\Delta'(1,1)
\end{equation} 
of motives over $\Q$.
\end{thm}

\medskip

Following \cite{ono}, we call $Y_\alpha$ {\it modular} if $L(h^2_\tr(Y_\alpha),s)$ is
the $L$-function of a Hecke eigenform of weight $3$.
The authors give a list of $\alpha$'s for $Y_\alpha$ to be modular.
\begin{thm}[{\cite[Theorem 1.2]{ono}}]\label{d=2-thm2}
The K3 surface $Y_\alpha$ is modular if and only if
$\alpha=-1,4^{\pm1},-8^{\pm1},64^{\pm1}.$
\end{thm}
The corresponding Hecke eigenforms are as follows (\cite[p.366--367]{ono}).
Let $\eta(z)$ be the Dedekind eta function. Let
\begin{equation}\label{d=2-table1}
A=\eta^6(4z),\quad B=\eta^2(z)\eta(2z)\eta(4z)\eta^2(8z),\quad
C=\eta^3(2z)\eta^3(6z),\quad D=\eta^3(z)\eta^3(7z)
\end{equation}
be weight $3$ newforms of level $16$, $8$, $12$, $7$ respectively.
Let $\chi_d$ denote the quadratic character associated to the quadratic field $\Q(\sqrt{d})$.
Then the corresponding Hecke eigenforms are given as follows.
\begin{equation}\label{d=2-table2}
\text{
 \begin{tabular}{c|ccccccc}
$\alpha$&$-1$&$4$&$1/4$&$-8$&$-1/8$&$64$&$1/64$\\
\hline
Hecke eigenform&$B\ot\chi_{-4}$&$C$&$C\ot\chi_{-4}$&$A$&$A\ot\chi_8$&$D$&$D\ot\chi_{-4}$
\end{tabular}
}
\end{equation}
where $f\ot\chi$ denotes the $\chi$-twist of the modular form.

\medskip

We turn to the K3 surface $X_\alpha$.
We discuss the Beilinson conjecture for
the real regulator map
\begin{equation}\label{d=2-eq1}
\reg_\R:K_3(X_\alpha)^{(3)}\lra 
H^2_B(X_{\alpha,\R},\R(2))^{F_\infty=1}
\end{equation}
where $X_{\alpha,\R}:=X_\alpha\times_\Q\R$.
By \eqref{XY},
\[
L(h^2_\tr(X_\alpha),s)=L(h^2_\tr(Y_\alpha),s).
\]
If $\alpha=-1,4^{\pm1},-8^{\pm1},64^{\pm1}$, then this is the $L$-function of
a Hecke eigenform by Theorem \ref{d=2-thm2}.

Let $f:X_\alpha\to\P^1(x_2)$ be the projection which gives an elliptic fibration.
Let $N(X_\alpha)\subset
\NS(X_\alpha)\ot\Q$ be the $\Q$-linear subspace generated
by irreducible components of the singular fibers and 
a section defined by $x_0=1$ and $x_1=\infty$. 
It is an easy exercise to show that the rank of $N(X_\alpha)$ is $19$.
Let $H^{pq}$ denote the Hodge $(p,q)$-component of 
$H^2_B(X_\alpha,\C)/N(X_\alpha)$.
Then $\dim H^{20}=\dim H^{20}=\dim H^{11}=1$ and
\begin{equation}\label{d=2-eq4}
\dim[H^2_B(X_\alpha,\Q)/N(X_\alpha)]^{F_\infty=\pm1}=
\overbrace{\dim(H^{20}+H^{02})^{F_\infty=\pm1}}^1
+\overbrace{\dim(H^{11})^{F_\infty=\pm1}}^{0\text{ or }1}.
\end{equation}
\begin{lem}\label{MW-lem}
For $\alpha\in \R\setminus\{0,1\}$, 
$\dim(H^{11})^{F_\infty=1}=0$ if and only if $\alpha>1$.
\end{lem}
\begin{pf}
It follows from Theorem \ref{GT} that there is an isomorphism
\[
H^2(X_\alpha,\Q)/\NS(X_\alpha)\cong H^1_B(E_\alpha,\Q)\ot H_B^1(E'_\alpha,\Q)/\Q[\Delta']
\]
of $3$-dimensional $F_\infty$-modules.
In general, $\dim H^1_B(X,\Q)^{F_\infty=1}
=\dim H^1_B(X,\Q)^{F_\infty=-1}$ for a projective smooth variety $X$ over $\R$.
From this, we have
\[
\dim[H^1_B(E_\alpha,\Q)\ot H^1_B(E'_\alpha,\Q)]^{F_\infty=+1}=
\dim[H^1_B(E_\alpha,\Q)\ot H^1_B(E'_\alpha,\Q)]^{F_\infty=-1}=2,
\] and hence
$\dim(H^{11})^{F_\infty=1}=0$ if and only if $\Q[\Delta']=(\Q[\Delta'])^{F_\infty=1}$,
which is equivalent to that 
the homomorphism $\delta':(x,y)\mapsto(x,y/\sqrt{1-\alpha})$
satisfies $F_\infty(\delta')=-\delta'$.
The last condition is equivalent to $1-\alpha<0$.
\end{pf}
Recall from Theorem \ref{per-thm1} the homology cycle $\Delta_\alpha$ with 
$0<|\alpha|\ll1$.
Extend $\Delta_\alpha$ for $\alpha\in\R\setminus\{0,1\}$ (this is not uniquely determined),
and put
\[
\Delta_\alpha^\pm:=\frac12(F_\infty\pm1)\Delta_\alpha\in 
H_2(U_\alpha,\Q)^{F_\infty=\pm1}.
\]
It follows from Theorem \ref{per-thm1}
that one has
\begin{equation}\label{d=2-eq3}
\int_{\Delta^+_\alpha} \omega_{1,1,1}=
\mathrm{Re}\left[{}_3F_2\left({\frac12,\frac12,\frac12\atop1,1};\alpha\right)\right]
\end{equation}
and this does not vanish 
\footnote{
cf. http://functions.wolfram.com/HypergeometricFunctions/Hypergeometric3F2/.}
for all $\alpha\in\R$.
This shows $\Delta^+_\alpha\ne0$.
By the construction of $\Delta_\alpha$,
the pairing $\langle -\mid\Delta_\alpha\rangle$
annihilates the subspace $\langle \fib,\infty\rangle$.
Hence we have a well-defined surjective homomorphism
\[
\langle-\mid\Delta^+_\alpha\rangle:
[H^2_B(X_\alpha,\Q)/N(X_\alpha)]^{F_\infty=1}\lra\Q,
\]
and this is bijective if and only if $\dim(H^{11})^{F_\infty=1}=0$.

\medskip

Summing up the above observation, we arrive at the following statement.
\begin{conj}[Weak Beilinson conjecture]
\label{d=2-WBconj}
Suppose $\dim(H^{11})^{F_\infty=1}=0$ or equivalently $\alpha>1$ by 
Lemma \ref{MW-lem}.
There is an integral element $\xi\in K_3(X_\alpha)^{(3)}_\Z$ such that
\begin{equation}\label{d=2-WBconj-eq1}
\frac{1}{(2\pi \sqrt{-1})^2}\langle\reg_\R(\xi)\mid\Delta^+_\alpha\rangle
\sim_{\Q^\times}L'(h^2_\tr(X_\alpha),0).
\end{equation}
Note $\langle\reg_\R(\xi)\mid\Delta^+_\alpha\rangle=\mathrm{Re}\langle\reg_\R(\xi)\mid\Delta_\alpha\rangle$.
\end{conj}
\begin{rem}
A K3 surface $X$ over a field of characteristic $0$
is called singular if $\dim\NS(X)=20$.
For such a $X$ over a number field, one can expect $K_i(X)^{(j)}
=K_i(X)^{(j)}_\Z$ for $i\ne1$.
Indeed there is the Shioda-Inose structure (\cite{I-S})
\begin{equation}\label{Shioda-Inose}
\xymatrix{
X_{\ol\Q}\ar[rd]_{\rho_Z}&&E\times E'\ar[ld]^{\rho_{E\times E'}}\\
&Z
}
\end{equation}
over $\ol\Q$,
in which the arrows are rational dominant maps of degree 2, and $Z$ is a K3 surface
and $E,E'$ are elliptic curves with complex multiplication.
Moreover $E$ and $E'$ are isogenous. Let $\xi\in K_i(X)^{(j)}$.
Then $\xi\in K_i(X)^{(j)}_\Z$ if and only if $\rho_{Z*}\xi\in K_i(Z)^{(j)}_\Z$, and
$\rho_{Z*}\xi\in K_i(Z)^{(j)}_\Z$ if and only if
$\rho_{E\times E'}^*\rho_{Z*}\xi\in K_i(E\times E')^{(j)}_\Z$.
Repalce the base field with a number field $F$ and we may assume that 
$E$ and $E'$ have everywhere good reduction.
Then $K_i(E\times E')^{(j)}_\Z$ is the kernel of the boundary map
\[
\bigoplus_\wp\partial_\wp:
K_i(E\times E')^{(j)}\lra\bigoplus_{\wp}K_{i-1}(E_\wp\times E'_\wp)^{(j-1)}
\]
where $\wp$ runs over all primes of $F$, and $E_\wp,E'_\wp$ denote the reduction at $\wp$.
However if $i\ne1$, a general conjecture says that the right hand side vanishies
(\cite[12.2]{janssen}).
Hence we have
$K_i(E\times E')^{(j)}=K_i(E\times E')^{(j)}_\Z$ for $i\ne1$ and the same statement for $X$.
\end{rem}
We examine the above conjecture for the higher Ross symbol
\[
\xi=\Ross=\left\{\frac{1-x_0}{1+x_0},\frac{1-x_1}{1+x_1},\frac{1-x_2}{1+x_2}\right\}.
\]

\begin{exmp}[$\alpha=4$]
By Theorem \ref{d=2-thm2} and \eqref{d=2-table2}, we have
$L(h^2_\tr(X_4),s)=L(C,s)$.
Apply Theorem \ref{mainB-thm}. We then have
\begin{align*}
\frac{1}{(2\pi \sqrt{-1})^2}\langle\reg_\R(\Ross|_{X_4})\mid\Delta^+_4\rangle&=
\mathrm{Re}[\cF_{\frac12,\frac12,\frac12}(t)|_{t=4}]\\
&=\mathrm{Re}\left[
-\log64+\log t
+\frac{t}{8}{}_5F_4\left({\frac32,\frac32,\frac32,1,1\atop2,2,2,2};t\right)\bigg|_{t=4}
\right].
\end{align*}
The computer calculation shows
\begin{align*}
&\mathrm{Re}[\cF_{\frac12,\frac12,\frac12}(t)|_{t=4}]= -2.41291989930352597175242344918,\\
&L'(C,0)=0.30161498741294074646905293114776839989
\end{align*}
and hence 
\[
\frac{1}{(2\pi \sqrt{-1})^2}\langle\reg_\R(\Ross|_{X_4})\mid\Delta^+_4\rangle=
-8L'(C,0)
\]
approximately.
\end{exmp}
\begin{exmp}[$\alpha=64$]
We have $L(h^2_\tr(X_{64}),s)=L(D,s)$
by Theorem \ref{d=2-thm2} and \eqref{d=2-table2}.
By Theorem \ref{mainB-thm},
\[
\frac{1}{(2\pi \sqrt{-1})^2}\langle\reg_\R(\Ross|_{X_{64}})\mid\Delta^+_{64}\rangle=
\mathrm{Re}[\cF_{\frac12,\frac12,\frac12}(t)|_{t=64}].
\]
The computer calculation shows
\begin{align*}
&\mathrm{Re}[\cF_{\frac12,\frac12,\frac12}(t)|_{t=64}]= -0.821372862231216089683652759186\\
&L'(D,0)=0.10267160777890201121045659489829291400.
\end{align*}
We have approximately
\[
\frac{1}{(2\pi \sqrt{-1})^2}\langle\reg_\R(\Ross|_{X_{64}})\mid\Delta^+_{64}\rangle=
-8L'(D,0).
\]
\end{exmp}

\subsection{Addendum : $X_\alpha$ for $\alpha=1$}\label{add-sect}
When $\alpha=1$ the surface $X_\alpha$ has a singular point,
but its smooth model is a K3 surface.
To discuss the weak Beilinson conjecture for $X_1$, we cannot
apply Theorem \ref{mainB-thm} directly, we need an additional argument.

\medskip

Let $U_t\to \Spec \Q[t,t^{-1}]$ be the hypergeometric scheme.
Take the base change by $t-1=\l^2$, and 
put $U:=\Spec \Q[[\l]][x_0,x_1,x_2]/((1-x_0^2)(1-x_1^2)(1-x_2^2)-1-\l^2)$.
Let $X^*\supset U$ be the projective scheme \eqref{sc-eq1} over $\Q[[\l]]$.
Let $\rho:X\to X^*$ be the blow-ups in Proposition \ref{sc-thm}.
Then $X\to\Spec\Q[[\l]]$ is smooth outside a point
$P=\{(x_0,x_1,x_2,\l)=(0,0,0,0)\}$. A neighborhood of $P$ in $X$ is locally defined by
an equation $x_0^2+x_1^2+x_2^2+\l^2=0$.
Let $\rho':\cX\to X$ be the blow-up at $P$, and we have a diagram
\[
\xymatrix{
\cX\ar[r]^{\rho \rho'}\ar[rd]_\pi&X^*\ar[d]\\
&\Spec\Q[[\l]].
}
\]
Then $\pi$ is a semistable family.
The central fiber of $\pi$ is $S\cup Q$ where $Q$ is a smooth
quadratic surface in $\P^3$, $L=Q\cap S$ is a smooth quadratic curve in $\P^2$
and $S$ is a K3 surface which is the smooth compactification of
an affine equation 
\[
(1-x_0^2)(1-x_1^2)(1-x_2^2)=1.
\]
Put $\cU:=(\rho\rho')^{-1}(U)$ and $U_S:=S\cap \cU$ and $i:U_S\hra \cU$ the inclusion.
Write $\Ross|_\cU:=(\rho\rho')^*(\Ross)\in K_3(\cU)$
the inverse image of the higher Ross symbol.
The restriction $i^*(\Ross|_\cU)$ lies in the image of $K_3(S)$, which we write by
$\Ross|_S$ simply.
Let $\Delta=\{|\l|<1\}$ be the unit disk, and $j:\Delta\setminus\{0\}\hra\Delta$
the inclusion.
The local monodromy $T$ on $H^2_B(\pi^{-1}(\l),\Q)$ is trivial since
so does on ${}_3F_2\left({\frac12,\frac12,\frac12\atop1,1};1-\l^2\right)$
(or this follows from \cite[p.118, algebraic monodromy criteria]{morrison}).
There is a commutative diagram
\[
\xymatrix{
K_3(\cX)\ar[r]\ar[d]_{i^*}&\Ext^1_{\VMHS(\Delta)}(\Q,j_*j^*R^2\pi_*\Q(3))\ar[d]^{i^*}
\ar[r]&\vg(\Delta\setminus\{0\},j^*R^2\pi_*\C/\Q(3))\ar[d]^{i^*}
\\
K_3(S)\ar[r]&\Ext^1_{\MHS}(\Q,H^2(S,\Q(3)))
\ar[r]&H^2(S,\C/\Q(3))
}
\]
where the right and middle vertical arrows are 
defined by the fact that the local monodromy is trivial.
In this situation, it is possible to apply Theorem \ref{mainB-thm} so that one has
\begin{equation}\label{add-eq0}
\frac{1}{(2\pi\sqrt{-1})^2}\langle\reg(\Ross|_S)\mid\Delta_1\rangle
=\cF_{\frac12,\frac12,\frac12}(t)|_{t=1}
\end{equation}
modulo $2\pi\sqrt{-1}\Q$ where
$\Delta_1\in H_2(U_S,\Q)$ is the homology cycle $\Delta_\alpha$ at $\alpha=1$.
Let $f:S\to \P^1(x_2)$ be the elliptic fibration.
Then the Neron-Severi group $\NS(S)\ot\Q$ is of rank $20$ generated by 
fibral divisors and a section defined by $(x_0,x_1)=(1,\infty)$.
The pairing 
\[
\langle-\mid\Delta^+_1\rangle:
[H^2_B(S,\Q)/\NS(S)\ot\Q]^{F_\infty=1}\os{\cong}{\lra} \Q
\]
with $\Delta_1^+$ is well-defined and bijective.
From \eqref{add-eq0} we have
\begin{align}
\frac{1}{(2\pi\sqrt{-1})^2}\langle\reg_\R(\Ross|_S)\mid\Delta^+_1\rangle
&=\mathrm{Re}[\cF_{\frac12,\frac12,\frac12}(t)|_{t=1}]\notag\\
&=-\log64
+\frac{1}{8}{}_5F_4\left({\frac32,\frac32,\frac32,1,1\atop2,2,2,2};1\right).\label{add-eq1}
\end{align}
Next we show $L(h^2_\tr(S),s)=L(A,s)$
where $A$ is the Hecke eigenform in \eqref{d=2-table1}\footnote{This is stated in \cite[p.354]{ono}, though
the proof is omitted.}.
Let $Z$ be the K3 surface defined by an equation
\[
(y_0-y_0^{-1})(y_1-y_1^{-1})(y_2+y_2^{-1})=8.
\]
Let $\rho_1:Z\to S$ be a covering given by
\[
\rho_1^*(x_0)=\frac12(y_0+y_0^{-1}),\quad
\rho_1^*(x_1)=\frac12(y_1+y_1^{-1}),\quad
\rho_1^*(x_2)=\frac{\sqrt{-1}}{2}(y_2-y_2^{-1}).
\]
Let $E_i$ be the elliptic curve over $\Q$ defined by an equation $w_i^2+z_i^4=1$.
Let $\rho_2:E_1\times E_2\to Z$ be a covering given by
\[
\rho_2^*(y_0)=\frac{z_2w_1}{z_1^2-1},\quad
\rho_2^*(y_1)=\frac{-w_1w_2+2z_1z_2}{w_1w_2+2z_1z_2},\quad
\rho_2^*(y_2)=\frac{z_1w_2}{z_2^2+1}.
\]
The composition $\rho:=\rho_1\circ\rho_2:E_1\times E_2\to S$ gives the Shoida-Inose structure over $\Q(\sqrt{-1})$. 
Let $\alpha_p,\beta_p$ be the eigenvalues of the $p$-th Frobenius $\Phi$ on 
the crystalline cohomology group $H^1_\crys(E_{i,p}/\Z_p)$ of the reduction
$E_{i,p}/\F_p$ at $p>2$ (note $E_1\cong E_2$).
Since
\[
\rho^*\left(\frac{dx_0dx_1}{(1-x_0^2)(1-x_1^2)x_2}\right)
=4\sqrt{-1}\frac{dz_1}{w_1}\frac{dz_2}{w_2},
\]
we have
\[\det(1-\Phi T\mid H^2_{\crys}(S_p/\Z_p))=
\begin{cases}
(1-\alpha_p^2T)(1-\beta_p^2T)&p\equiv 1\text{ mod }4\\
1+p^2T^2&p\equiv 3\text{ mod }4.
\end{cases}
\]
This implies
\begin{equation}\label{add-eq2}
L(h^2_\tr(S),s)=L(h^2_\tr(E_1\times E_2),s)=L(A,s).
\end{equation}
Finally we employ a formula of D. Samart \cite[Corollary 1.3]{samart1}
\begin{equation}\label{add-eq3}
\log64
-\frac{1}{8}{}_5F_4\left({\frac32,\frac32,\frac32,1,1\atop2,2,2,2};1\right)
=8L'(A,0).
\end{equation}
By \eqref{add-eq1}, \eqref{add-eq2} and \eqref{add-eq3},
we deduce the weak Beilinson conjecture for $\Ross|_S$.
\begin{thm}\label{samart-thm}
\[
\frac{1}{(2\pi\sqrt{-1})^2}\langle\reg_\R(\Ross|_S)\mid\Delta^+_1\rangle
=-8L'(h^2_\tr(S),0).
\]
\end{thm}

\end{document}